\documentclass[10pt,letterpaper]{article}
\usepackage[top=1in,left=1.2in,footskip=1in,marginparwidth=1in]{geometry}

\usepackage{multirow}
\usepackage[utf8]{inputenc}
\usepackage{color, colortbl}
\usepackage{amsmath}
\usepackage{amssymb}
\usepackage{algorithm}
\usepackage{float}
\usepackage[noend]{algpseudocode}
\usepackage[misc]{ifsym}
\newcolumntype{a}{>{\columncolor{gray!20}}c}
\usepackage{caption}
\usepackage{xcolor}
\usepackage{color}

\usepackage{ bbold }

\usepackage{mleftright}

\newcommand{\real}{\mathbb{R}}

\usepackage{tikz}
\usetikzlibrary{matrix, fit}
\usetikzlibrary{backgrounds}

\usepackage{cite}

\usepackage{nameref,hyperref}

\usepackage[right]{lineno}

\usepackage{microtype}

\usepackage{changepage}

\usepackage[aboveskip=1pt,labelfont=bf,labelsep=period,singlelinecheck=off]{caption}

\makeatletter
\renewcommand{\@biblabel}[1]{\quad#1.}
\makeatother

\usepackage{color}

\definecolor{Gray}{gray}{.25}

\usepackage{graphicx}
\usepackage[colorinlistoftodos]{todonotes}
\usepackage{sidecap}

\usepackage{wrapfig}
\usepackage[pscoord]{eso-pic}
\usepackage[fulladjust]{marginnote}
\reversemarginpar

\begin{document}
\vspace*{0.35in}

\begin{flushleft}
{\Large
\textbf\newline{Surface Eigenvalues with Lattice-Based Approximation}}\\
%
%
In comparison with analytical solution
\newline
\\
Y.~Wu\textsuperscript{1,\Letter},
T.~Wu\textsuperscript{1},
and S.-T.~Yau\textsuperscript{1,2}
\\
\bigskip
$^1$Center of Mathematical Sciences and Applications, Harvard University\\
$^2$Department of Mathematics, Harvard University\\
$^{\text{\Letter}}${e-mail: ywu@cmsa.fas.harvard.edu}

\end{flushleft}


\tableofcontents

\begin{abstract}
In this paper, we propose a meshless method of computing eigenvalues and eigenfunctions of a given surface embedded in $\mathbb R^3$. We use point cloud data as input and generate the lattice approximation for some neighborhood of the surface. We compute the eigenvalues and eigenvectors of the cubic lattice graph as an approximation of the eigenvalues and eigenfunctions of the Laplace-Beltrami operator on the surface. We perform extensive numerical experiments on surfaces with various topology and compare our computed eigenvalues from point cloud surface with exact solutions and standard finite element methods using triangle mesh.
\end{abstract}

\section{Introduction}
For a given surface, the associated Laplace-Beltrami operator is a useful tool that describes in part the intrinsic geometric information  of the shape. 
In particular, the eigenvalues and eigenfunctions of the Laplace-Beltrami operator, i.e., the spectrum information, are employed as global shape descriptors.

Computing eigenvalues and eigenfunctions for surfaces has a wide range of applications in mathematics and computer graphics.
In mathematics, it has long been a fundamental problem to study the connections between eigenvalues and other geometric quantities such as genus and the Cheeger constants \cite{yau2000review}. 
Since the analytic solutions of most eigenvalue problems are only known for special cases and usually involve meticulous analysis solving partial differential equations, a numerical algorithm for computing the eigenvalues and eigenfunctions could play an important role in discovering and verifying mathematical theorems on eigenvalue problems.

In computer graphics, computing eigenvalues and eigenfunctions is a fundamental step in spectral shape analysis. The eigenvalues and eigenfunctions could be used to measure the difference between the intrinsic geometries of the two surfaces. This is particularly useful for non-rigid shapes in 3-dimensional space.
Furthermore, eigenvalues and eigenvectors could be used for surface registrations \cite{hamidian2019surface}. See \cite{belkin2003laplacian} \cite{coifman2006diffusion} \cite{reuter2009discrete}  \cite{reuter2005laplace}\cite{seo2010heat} for examples of such applications.

\subsection{Contributions}

Most existing methods for computing eigenvalues and eigenfunctions rely on the triangle mesh approximation of a surface; see \cite{reuter2006laplace}\cite{hamidian2019surface}. 
The advantages of the point cloud-based eigenvalue method are:

\begin{itemize}
\item Point cloud data uses only coordinate information, which can be regarded as a subset of the information contained in a triangle mesh as the set of vertices. Comparing the accuracy of computed eigenvalues between mesh-based algorithm and meshless algorithm reveals in part the impact of surface information, which is of great theoretical importance.

\item Collecting point cloud data is convenient. Contemporary 3D scanners provide 3D point cloud data sampled from the surfaces of solid objects, and mobile devices with depth sensors can collect point cloud data also. Since generating meshes from point cloud may not be feasible, algorithms computing eigenvalues directly from point cloud data are useful.

\item The concept of lattice approximation could be easily generalized to higher-dimensional cases and less regular geometric objects, such as orbifolds.


\end{itemize}


 In this article, we approximate a surface $M\subset\mathbb R^3$ by a point cloud, and we compute the eigenvalues of the Laplacian on a surface $M$ by approximating the $\epsilon$-neighborhood of the surface, then we compute the discrete Laplacian using the lattice. 
We visualize the convergence of computed discrete eigenvalues toward exact eigenvalues solved analytically. 

\subsection{Related works}

Reuter et al. \cite{reuter2006laplace} introduced the spectrum of the Laplace-Beltrami operator of a shape as the signature or fingerprint of the shape. Since the spectrum is an isometry invariant, it is independent of the object’s representation, including parametrization and spatial position.
The authors pointed out that in spectrum methods, checking if two objects are isometric needs no prior alignment (registration/localization) of the objects but only a comparison of their spectra. Moreover, the authors demonstrated that it is computationally feasible to extract elementary geometrical data such as the volume, the boundary length, and even the Euler characteristic from the numerically calculated eigenvalues.  This indicates the geometrical importance of the eigenvalues and eigenfunctions. It is also suggested in \cite{reuter2006laplace} to implement spectrum methods to support copyright protection, database retrieval, and quality assessment of digital data that represent surfaces and solids. Rustamov \cite{rustamov2007laplace} proposed to use the eigenvalues and eigenfunctions of surfaces to do shape clustering and classification. Le\'vy \cite{levy2006laplace} studied the specific type of function bases defined by the eigenfunctions of the Laplace-Beltrami operator. Such a  function basis is well-adapted to the geometry and the topology of the object.

Dong et al. \cite{dong2005quadrangulating}, 
Le\'vy \cite{levy2006laplace}, and
Reuter \cite{reuter2010hierarchical} used eigenfunctions to develop approaches to segment shapes.
The resulting patches are well-shaped and arise naturally from the intrinsic properties of the surface. Such approaches work for any surface, regardless of its genus. In \cite{reuter2010hierarchical}, algorithms are developed to segment surfaces into meaningful parts and to register these parts across populations of near intrinsic isometric shapes, such as heads, arms, legs, and fingers of humans in different body postures. 
The method utilized the fact that quantitative and qualitative behaviors of the eigenvalues and eigenfunctions are similar for near intrinsic isometric shapes.  Topological features (level sets and Morse theory) are implemented for the segmentation. For the purpose of  eliminating topological noise and comparing eigenfunctions, concepts from persistent homology are employed. \cite{reuter2010hierarchical} also discussed the computation of eigenfunctions and eigenvalues using cubic finite elements on triangle meshes and  constructing persistence diagrams by the Morse-Smale complex.

In \cite{reuter2010hierarchical}, the registration of the shapes is mainly for near isometric surfaces. In reality, it is often needed to find the registration of surfaces that are not near isometric. For example, heart motion and brain development are not isometric. Therefore, it could be challenging to implement the Laplace spectrum methods in these applications.
Shi et al. \cite{shi2011conformal} developed a novel technique for surface deformation and mapping in the high-dimensional Laplace-Beltrami embedding space.
For a surface, the authors deformed its Laplace-Beltrami eigenvalues and eigenfunctions and
realized its deformation in the high-dimensional Laplace-Beltrami embedding space by iteratively optimizing conformal metrics. By this deformation technique, 
the authors developed an approach for surface mapping between non-isometric surfaces
and demonstrated its application in mapping hippocampal atrophy of multiple sclerosis (MS) patients with depression \cite{shi2011conformal}. 
Hamidian et al. \cite{hamidian2016quantifying}\cite{hu2016visualizing} also described an approach using eigenvalues and eigenfunctions to quantify and visualize non-isometric deformations of surfaces. Such deformation could be expressed as a linear interpolation of eigenvalues of the two surfaces, realized by a time-dependent scale function defined on each vertex. Each iteration amounts to solving a quadratic programming problem based on a spectrum variation theorem and a smoothness energy constraint. The final scale function can be obtained by combining the deformations from each step and quantitatively describing non-isometric deformations between two shapes. 

There are other works that employ the spectrum of Laplacian for surface matching. Rodol\`a et al. \cite{rodola2017partial} 
used perturbation analysis to show how the Laplace-Beltrami eigenfunctions change after removing part of the shape and exploited it as a prior on the spectral representation of the correspondence between non-rigid shapes.
Litany et al. \cite{litany2017fully} extended the study of partial matching, including the presence of topological noise.
Kovnatsky et al. \cite{kovnatsky2013coupled} showed how to modify (align) the eigenvectors of the Laplace-Beltrami operator to match non-isometric surfaces. Ovsjanikov et al. \cite{ovsjanikov2012functional}
described a spectral method for shape matching by finding an alignment between eigenfunctions based on linear constraints. In \cite{rustamov2013map}\cite{zeng2010shape}, visualization of shape deformations based on spectral representations of the correspondence was shown.
In \cite{hamidian2019surface}, the authors
used both eigenvalues and eigenvectors to align two manifolds and extracted feature points from eigenvectors and employed them to align two surfaces.\\

\noindent \bf Acknowledgments. \rm The authors would like to thank Cliff Taubes for a great deal of inspiration from his work and fruitful discussions regarding the analytic solution and all aspects of this paper. The computations in this paper were run on the FASRC Cannon cluster supported by the FAS Division of Science Research Computing Group at Harvard University. The work is also supported in part by NSF 1760471.

\section{Analytic Solution of Surface Laplace’s Equation}

In this section, we present analytic solutions of Laplace’s equation of a sphere and a cone. Since the analytic solutions provide the exact values of eigenvalues of surface Laplacian, they can serve as the ground truth for the evaluation and comparison of the accuracy and convergence of numerical eigenvalues obtained by different algorithms.

\subsection{Spherical Harmonics}

Spherical harmonics are functions defined on the surface of a sphere. They form a complete set of orthogonal functions and thus an orthonormal basis. Each function defined on the surface of a sphere can be written as a sum of these spherical harmonics. Spherical harmonics may be organized by spatial angular frequency, and they are basis functions for irreducible representations of $SO(3)$.

Spherical harmonics can be defined as homogeneous polynomials of degree $\ell$ in $(x,y,z)$ that obey Laplace's equation.
A specific set of spherical harmonics, denoted ${\displaystyle Y_{\ell }^{m}(\theta ,\varphi )}$, forms an orthogonal system and is known as Laplace's spherical harmonics.
Laplace's equation imposes that the Laplacian of a scalar field $f:\mathbb {R} ^{3}\to \mathbb {C}$ is zero. In spherical coordinates, this is \cite{courant2008methods}:
$$\nabla ^{2}f={\frac {1}{r^{2}}}{\frac {\partial }{\partial r}}\left(r^{2}{\frac {\partial f}{\partial r}}\right)+{\frac {1}{r^{2}\sin \theta }}{\frac {\partial }{\partial \theta }}\left(\sin \theta {\frac {\partial f}{\partial \theta }}\right)+{\frac {1}{r^{2}\sin ^{2}\theta }}{\frac {\partial ^{2}f}{\partial \varphi ^{2}}}=0$$
with angular solutions:
$$Y_{\ell }^{m}:S^{2}\to \mathbb {C} $$
as a product of trigonometric functions.
To find solutions of the form $f(r, \theta, \varphi) = R(r) Y(\theta, \varphi)$, we use separation of variables resulting in two differential equations:
\begin{align*}
& {\frac {1}{R}}{\frac {d}{dr}}\left(r^{2}{\frac {dR}{dr}}\right)=\lambda, \\
& {\frac {1}{Y}}{\frac {1}{\sin \theta }}{\frac {\partial }{\partial \theta }}\left(\sin \theta {\frac {\partial Y}{\partial \theta }}\right)+{\frac {1}{Y}}{\frac {1}{\sin ^{2}\theta }}{\frac {\partial ^{2}Y}{\partial \varphi ^{2}}}=-\lambda.
\end{align*}
This system of equations provides a solution of the Laplace equation:
\begin{align*}
\nabla ^{2}f&={\frac {1}{r^{2}}}{\frac {\partial }{\partial r}}\left(r^{2}{\frac {\partial RY}{\partial r}}\right)+{\frac {1}{r^{2}\sin \theta }}{\frac {\partial }{\partial \theta }}\left(\sin \theta {\frac {\partial RY}{\partial \theta }}\right)+{\frac {1}{r^{2}\sin ^{2}\theta }}{\frac {\partial ^{2}RY}{\partial \varphi ^{2}}}=0.\end{align*}
$Y^m_\ell$ can be represented as a complex exponential, and associated Legendre polynomials:
$$Y_{\ell }^{m}(\theta ,\varphi )=Ne^{im\varphi }P_{\ell }^{m}(\cos {\theta }),$$
that fulfill:
$$\displaystyle r^{2}\nabla ^{2}Y_{\ell }^{m}(\theta ,\varphi )=-\ell (\ell +1)Y_{\ell }^{m}(\theta ,\varphi ).$$
Here $\displaystyle Y_{\ell }^{m}:S^{2}\to \mathbb {C}$ is called a spherical harmonic function of degree $\ell$ and order $m$, $P_{\ell }^{m}:[-1,1]\to \mathbb {R}$ is an associated Legendre polynomial, $N$ is a normalization constant, and $\theta \in [0,\pi]$ and $\varphi\in[0,2\pi)$ represent colatitude and longitude, respectively. For a fixed integer $\ell$, every solution $Y^\ell_m(\theta, \varphi) :S^{2}\to \mathbb {C}$, of the eigenvalue problem
\begin{equation}
r^{2}\nabla ^{2}Y^\ell_m=-\ell (\ell +1)Y^\ell_m
\end{equation}
is a linear combination of $Y^\ell_m: S^{2}\to \mathbb {C}$. There are $2\ell + 1$ linearly independent such polynomials. The multiplicity of the eigenvalues on a sphere equals the dimension of the space of homogenous, harmonic polynomials of degree $\ell$ \cite{shubin1987pseudodifferential}. Therefore, the multiplicity of the eigenvalue $-\ell(\ell+1)$ is $2\ell + 1$.

\subsection{Analytic Solution of Surface Laplace's Equation on a Cone}

\begin{figure}[!htb]
\begin{center}
\includegraphics[width=0.32\textwidth]{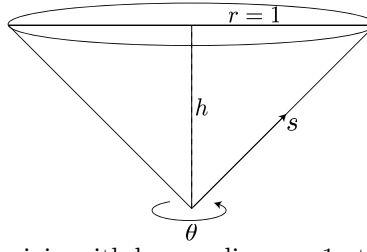}
\caption{Cone centered at the origin with base radius $r=1$ at height $z = 1$.}\label{fig:cone}
\end{center}
\end{figure}

In this section, we first present the analytic solution of Laplace's equation on the surface of a cone centered at the origin with base radius $r=1$ at height $z = 1$, as illustrated in Fig.~\ref{fig:cone}. 
We solve for the analytic solution of surface Laplace's equation on a cone by equating Laplace's equation and its first-order derivative on the side and top of the cone.

\noindent\textbf{Equation on the side of the cone.}\\

We parameterize the side of the cone as:
$$(t, \theta) \mapsto (t \cos\theta, t \sin \theta, ht), \quad 0 \leq t \leq 1, \quad 0 \leq \theta \leq 2\pi.$$
Then we have the metric matrix:
$$g 
= \begin{pmatrix}
dt & d\theta\\
\end{pmatrix}
\begin{pmatrix}
\langle g_{t},g_t\rangle & \langle g_{t},g_\theta\rangle\\
\langle g_\theta, g_t\rangle & \langle g_\theta,g_\theta\rangle\\
\end{pmatrix}
\begin{pmatrix}
dt\\
d\theta
\end{pmatrix}
= ( 1+h^2 )dt^2 + t^2 d^2\theta.$$

Now we change the variables:
$$s =t \sqrt{1 + h^2}  \Rightarrow 
g = ds^2 + {s^2 \over 1 + h^2} d^2\theta,$$
and the respective first fundamental form can be found in:
$$g = ds^2 + {s^2 \over 1 + h^2} d^2\theta
= \begin{pmatrix}
ds & d\theta\\
\end{pmatrix}
\begin{pmatrix}
1 & 0\\
0 &  {s^2 \over 1 + h^2}\\
\end{pmatrix}
\begin{pmatrix}
ds\\
d\theta
\end{pmatrix}.$$
The inverse metric tensor is:
$$g^{-1} = \begin{pmatrix}
1 & 0\\
0 &  {1 + h^2 \over s^2}\\
\end{pmatrix}.$$
The Laplacian is:
$$\Delta = {1\over \sqrt{\det g}} {\partial \over \partial x^i}\left(g^{ij} \sqrt{\det g}{\partial \over \partial x^i}\right),$$
so in our case,
\begin{align*}
    \Delta 
&= {1\over s} {\partial \over \partial s}+ {\partial^2 \over \partial s^2}
+ {1+h^2\over s^2} {\partial^2 \over \partial \theta^2}.
\end{align*}

Therefore, to find the Laplace eigenvalues on the side of a cone, we solve:
$${1\over s} f_s + f_{ss} + {1+h^2\over s^2}f_{\theta\theta} = -Ef.$$
Since $f$ is a smooth function, we can write it as:
$$f = \sum_{n =0} a_n(s)e^{in\theta}.$$
Since $e^{in\theta}$s are linearly independent, this is equivalent to solving:
$${1\over s} \left(a_n(s)e^{in\theta}\right)_s +  \left(a_n(s)e^{in\theta}\right)_{ss} + {1+h^2\over s^2} \left(a_n(s)e^{in\theta}\right)_{\theta\theta} = -E \left(a_n(s)e^{in\theta}\right).$$
So we have:
\begin{equation}\label{eqn:bessel}
{1\over s} (a_n)_s + (a_n)_{ss} - {n^2(1+h^2)\over s^2}(a_n)_{\theta\theta} = -Ea_n    
\end{equation}
for each $n$. 

The Bessel differential equation is the linear second-order ordinary differential equation given by:
$$x^2 {d^2y\over dx^2} + x {dy\over dx} + (x^2 - n^2)y = 0.$$
Therefore,
rearranging Equation (\ref{eqn:bessel}) gives the Bessel differential equation:
$$s (a_n)_s + s^2(a_n)_{ss} - \left(n^2(1+h^2) + Es^2\right)(a_n)_{\theta\theta} = 0$$
with solution $J_x\left(\sqrt{E}s\right)$ for $x^2 = (1+h^2)n^2 = 2n^2$.
Change of variable by $x = s\sqrt{E}$, so ${d\over dx} = {d\over d (s\sqrt{E})}={1\over \sqrt{E}} {d\over d s}$ and ${d^2\over dx^2}={1\over \sqrt{E}} {d\over d s}{1\over \sqrt{E}} {d\over d s}={1\over E} {d^2\over d s^2}$.
Then
$$s^2E{1\over E} {d^2 a_n\over d s^2}
+ s\sqrt{E} {1\over \sqrt{E}} {d a_n\over d s}
+ (s^2E - n^2)a_n = 0,$$
which simplifies to:
\begin{equation}\label{eqn:side}
s^2{d^2 a_n\over d s^2} + s {d a_n\over d s}+ (s^2E - n^2)a_n = 0.
\end{equation}

\vskip 3mm

\noindent
\textbf{Equation on the top of the cone.}\\

We parameterize the top of the cone as:
$$(r, \theta) \mapsto (r \cos\theta, r \sin \theta, 1), \quad 0 \leq r \leq 1, \quad 0 \leq \theta \leq 2\pi.$$
Then the metric matrix is:
$$g 
= \begin{pmatrix}
dr & d\theta\\
\end{pmatrix}
\begin{pmatrix}
\langle g_{r},g_r\rangle & \langle g_{r},g_\theta\rangle\\
\langle g_\theta, g_r\rangle & \langle g_\theta,g_\theta\rangle\\
\end{pmatrix}
\begin{pmatrix}
dr\\
d\theta
\end{pmatrix}
= \begin{pmatrix}
dr & d\theta\\
\end{pmatrix}
\begin{pmatrix}
1 & 0\\
0 & r^2\\
\end{pmatrix}
\begin{pmatrix}
dr\\
d\theta
\end{pmatrix}
= dr^2 + r^2 d^2\theta.$$

So the determinant of the first fundamental form is $r^2$, and the inverse metric tensor is:
$$g^{-1} = \begin{pmatrix}
1 & 0\\
0 &  {1 \over r^2}\\
\end{pmatrix}.$$
Thus, the Laplacian is:
\begin{align*}
    \Delta 
    &={1\over r} {\partial \over \partial r}\left(r{\partial \over \partial r}\right) + {1\over r} {\partial \over \partial \theta}\left({1\over r}{\partial \over \partial \theta}\right).
\end{align*}

Therefore, to find the Laplace eigenvalues on the top of a cone, we solve:
$${1\over r} v_r + v_{rr} + {1\over r^2}v_{\theta\theta} = -Ev.$$
Since $v$ is a smooth function, we can write it as:
$$v = \sum_{m =0} b_m(s)e^{im\theta}.$$
Since $e^{im\theta}$s are linearly independent, this is equivalent to solving:
\begin{equation}\label{eqn:solve}
{1\over r} \left(b_me^{im\theta}\right)_r + \left(b_me^{im\theta}\right)_{rr} + {1\over r^2}\left(b_me^{im\theta}\right)_{\theta\theta} = -E\left(b_me^{im\theta}\right).
\end{equation}
Rearranging gives the Bessel differential equation:
\begin{equation}\label{eqn:top}
    r (b_m)_r + r^2(b_m)_{rr} + \left( Er^2-m^2\right)(b_m)_{\theta\theta} = 0
\end{equation}
with a solution $J_m\left(\sqrt{E}r\right)$, so $b_m =\alpha J_m\left(\sqrt{E}r\right)$ for some constant $\alpha$.

Because the solution to the respective Laplace equations on the side of the cone and on the top of the cone have to agree on the edge to be continuous, the solutions to Equations (\ref{eqn:side}) and (\ref{eqn:top}) have to be the same at $r=1$. There is also a continuity condition on their derivatives at $r=1$, that is, $a_n = b_m$ for $m = n$.
\begin{align*}
    a_n &= \left.J_{\sqrt{2}n}\left(\sqrt{E}s\right)\right|_{s=t\sqrt{1+h^2} \text{ at } t=1}=J_{\sqrt{2}n}\left(\sqrt{2E}\right),\\
    b_m &=\left.\alpha J_m\left(\sqrt{E}r\right)\right|_{r=1}
    =\alpha J_m\left(\sqrt{E}\right).
\end{align*}
This gives the restraint:
\begin{equation}\label{eqn:sol_1}
    J_{\sqrt{2}n}\left(\sqrt{2E}\right)=\alpha J_n\left(\sqrt{E}\right).
\end{equation}
And also their first derivatives:
\begin{align*}
    a^\prime_n &=J^\prime_{\sqrt{2}n}\left(\sqrt{2E}\right),\\
    b^\prime_m &=\alpha J^\prime_m\left(\sqrt{E}\right),
\end{align*}
which gives the restraint:
\begin{equation}\label{eqn:sol_2}
   -J^\prime_{\sqrt{2}n}\left(\sqrt{2E}\right)=\alpha J^\prime_n\left(\sqrt{E}\right)
\end{equation}
because ${\partial\over\partial r}$ points outward from the boundary of the disk, so it corresponds to  $-{\partial\over\partial s}$, which points in from the top of the cone.
Bring Equation (\ref{eqn:sol_1}) into Equation (\ref{eqn:sol_2}) gives
the analytic solution of surface Laplacian on a cone:
\begin{equation}\label{eqn:solution}
J^\prime_{\sqrt{2}n}\left(\sqrt{2E}\right)J_n\left(\sqrt{E}\right)+ J_{\sqrt{2}n}\left(\sqrt{2E}\right)J^\prime_n\left(\sqrt{E}\right)=0.
\end{equation}
Using the recurrence relation:
\begin{equation}\label{eqn:rec}
J^\prime_v(z) = -J_{v+1}(z) + {v\over z}J_v(z),
\end{equation}
we obtain the analytic solution, displayed in Table \ref{tab:cone}.

\vskip 3mm

\noindent
\textbf{Multiplicities.}

\vskip 3mm

According to Equation (\ref{eqn:solve}), each solution $E$ of Equation (\ref{eqn:solution}) corresponds to two eigenvalues after expanding $e^{im\theta}$ into the real part and imaginary part for all $n>0$. Therefore, the solution from Bessel $J_0$ functions has multiplicity 1, and all others have multiplicity 2. Analytic solutions are exhibited in Table \ref{tab:bessel_eigenvalue}.

\section{Methods of Computing Eigenvalues of Surfaces}
In this section, we describe the classical method \cite{reuter2009discrete} for computing eigenvalues and eigenfunctions of surfaces embedded in $\mathbb R^3$ (Algorithm 1). Then we introduce an algorithm computing eigenvalues using a lattice that is approximated based on the triangle mesh (Algorithm 2), followed by our main algorithm, which computes eigenvalues using a lattice that is approximated based on the point cloud (Algorithm 3). A recent embedding method, the Closet Point Method for computing solutions to a variety of partial differential equations based on closest point representation of surface \cite{macdonald2011solving}, is analogous to our lattice approximation. Similarly, \cite{brandman2008level} solved eigenvalues of an elliptic operator defined on a compact hypersurface in $\real^n$ via computing an elliptic eigenvalue problem in a bounded domain on the Cartesian grid.

\subsection{Classical Method Using Triangle Mesh}

We first describe the classical method of computing eigenvalues on a triangle mesh surface by Reuter \cite{reuter2009discrete}. Given a triangle mesh $T=(V,E,F)$ embedded in $\mathbb R^3$, one of the most well-known notions of discrete Laplacian on this piecewise linear surface is the so-called cotangent Laplacian. For any $f\in\mathbb R^V$, its cotangent Laplacian $L$ is defined as
$$
(Lf)_i=\sum_{j:ij\in E}w_{ij}(f_j-f_i),
$$
where 
$$
w_{ij}=\cot\alpha_{ij}+\cot\beta_{ij}
$$
where $\alpha_{ij}$ and $\beta_{ij}$ are two inner angles in $(T,E,F)$ that are facing the edge $ij$.

\begin{algorithm}
\caption{Computing eigenvalues and eigenvectors of surfaces}\label{alg:triangle}

\textbf{Input:}  A triangle mesh $T=(V,E,F)$ and the number $k$ of wanted eigenvalues and eigenvectors.
\\
\textbf{Output:} First $k$ non-zero eigenvalues and associated eigenvectors.

\begin{algorithmic}[1]
\State Compute the inner angles of $T=(V,E,F)$.

\State Compute the (normalized) sparse matrix of the cotangent Laplacian $L$.

\State Compute the first $k$ non-zero eigenvalues and associated eigenvectors of $L$.
\end{algorithmic}
\end{algorithm}

\subsection{Lattice Approximation Method Using Triangle Mesh} 

Given a triangle mesh $T=(V,E,F)$ embedded in $\mathbb R^3$, for a relatively small $r>0$ and a sufficiently dense cubic lattice $(\mathbb Z/n)^3$, we can take $V^\prime = B_r(T)\cap(\mathbb Z/n)^3$ as a lattice approximation. Then, we approximate the surface Laplacian by the standard graph Laplacian on the standard cubic lattice graph upon $V^\prime$ defined as:
\begin{equation}
    \label{mask_function}
    L_{i,j} = 
    \begin{cases}
        \deg(v_1) &  \text{ if } i = j\\
        -1 &\text{ if } i \neq j \text{ and } v_i \text{ is adjacent to } v_j\\
        0  & \text{ otherwise.}\\
    \end{cases},
\end{equation}

\begin{algorithm}
\caption{Computing eigenvalues and eigenvectors of surfaces}\label{alg:triangle}

\textbf{Input:}  A triangle mesh $T=(V,E,F)$, and a radius parameter $r>0$, and a lattice density parameter $n\in\mathbb Z^+$, and the number $k$ of wanted eigenvalues and eigenvectors.
\\
\textbf{Output:} First $k$ non-zero eigenvalues and associated eigenvectors.

\begin{algorithmic}[1]
\State Compute the set $V^\prime$ of lattice vertices in the $r$-neighborhood of the triangle mesh $T$.

\State Compute the (normalized) sparse matrix of the graph Laplace $L$.

\State Compute the first $k$ non-zero eigenvalues and associated eigenvectors of $L$.
\end{algorithmic}
\end{algorithm}

\newpage

\subsection{Meshless Algorithm}

Given a point cloud  $P$ in $\mathbb R^3$, then for a relatively small $r>0$ and a sufficiently dense cubic lattice $(\mathbb Z/n)^3$, we can take $V = B_r(P)\cap(\mathbb Z/n)^3$ as a lattice approximation. Further, we approximate the surface Laplacian by the standard graph Laplacian on the standard cubic lattice graph upon $V$.

\begin{algorithm}[!htb]
\caption{Computing eigenvalues and eigenvectors of surfaces}\label{alg:meshless}

\textbf{Input:}  A point cloud $P$, and a radius parameter $r>0$, and a lattice density parameter $n\in\mathbb Z^+$, and the number $k$ of wanted eigenvalues and eigenvectors.
\\
\textbf{Output:} First $k$ non-zero eigenvalues and associated eigenvectors.

\begin{algorithmic}[1]
\State Compute the set $V$ of lattice vertices as described in Section 2.1.

\State Compute the (normalized) sparse matrix of the graph Laplace $L$.

\State Compute the first $k$ non-zero eigenvalues and associated eigenvectors of $L$.
\end{algorithmic}
\end{algorithm}

\section{Numerical Results on Surface Eigenvalues}
\subsection{Convergence About $r$-neighborhood and Lattice Density on Sphere}
   
\begin{figure}[!htb]
\begin{center}
\includegraphics[width=0.36\textwidth]{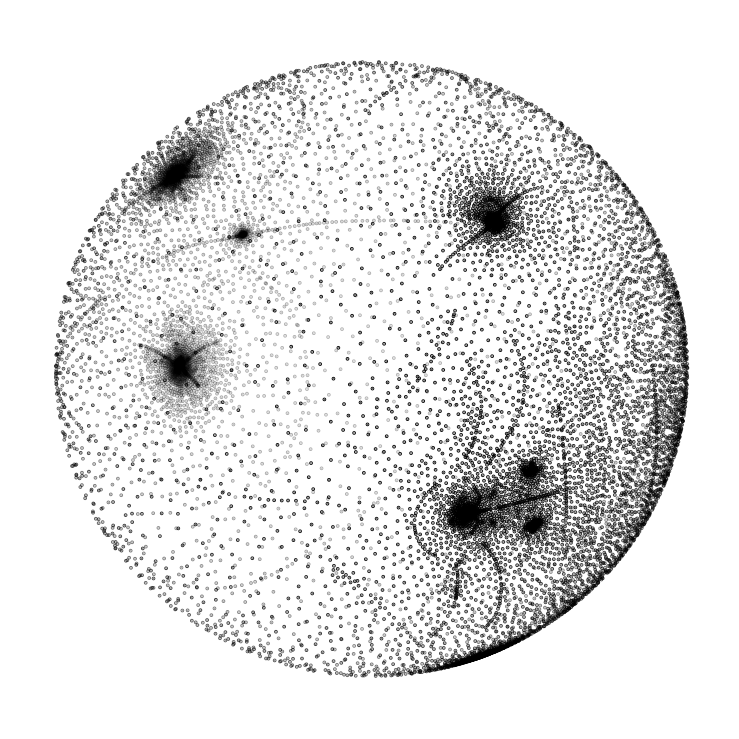}\hskip 15mm\includegraphics[width=0.342\textwidth]{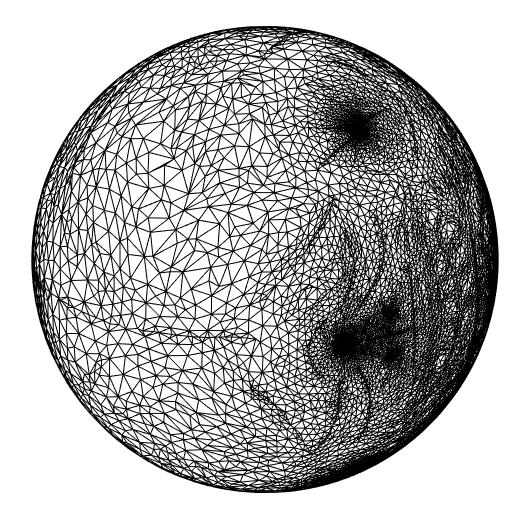}
\caption{(\it Left\rm) point cloud of $S^2$; (\it right\rm) triangle mesh of $S^2$.}\label{fig:s2}
\end{center}
\end{figure}

We display the computed first ten non-zero eigenvalues of $S^2$. Further experiments indicate that Algorithm \ref{alg:triangle} and Algorithm \ref{alg:meshless} give good approximations of the eigenvalues of $S^2$, as $n\rightarrow\infty$ for properly chosen $r>0$. We also did numerical experiments for different point cloud approximations of $S^2$, and converging speeds were different for different point clouds. The point cloud we used for $S^2$ is exhibited in Fig.~\ref{fig:s2}.

Table \ref{tab:eig_alg_triangle_nr} reports the first ten non-zero surface eigenvalues on $S^2$ (counting multiplicity) by Algorithm \ref{alg:triangle} over grid size ${1\over n}$ for $n=40$ to 110 over shell thickness $r = 0.02 $ to 0.05 over $S^2$,
which is the $\epsilon$-neighborhood of the grid. 
 According to Table \ref{tab:eig_alg_triangle_nr}, we see that as the shell becomes thinner (smaller $y$ value) and the grid becomes finer (bigger $x$ value), the eigenvalues of the discrete $S^2$ converge to the true eigenvalues on a sphere, as discussed in Section 2.1. The multiplicity of eigenvalues $-\ell(\ell+1)$ is $2\ell + 1$ for $\ell = 1,2,3$, as expected.
  
Table \ref{tab:eig_alg_meshless_nr} reports the first ten non-zero surface eigenvalues by Algorithm \ref{alg:meshless} over grid size ${1\over n}$ for $n=40$ to 110 over shell thickness $r = 0.05 $ to 0.08 over $S^2$. Similar trends are found between 
Table \ref{tab:eig_alg_triangle_nr} and Table \ref{tab:eig_alg_meshless_nr}. The
results suggest that as the thickness of a surface and the length of the cubes approach zero, the computed discrete eigenvalues approximate the smooth eigenvalue more accurately.

Fig.~\ref{fig:ev_triangle_3D} and Fig.~\ref{fig:ev_meshless_3D} exhibit computed surface eigenvalues by Algorithm \ref{alg:triangle} and Algorithm \ref{alg:meshless}, respectively.
Fig.~\ref{fig:ev_meshless_radius} displays the change of eigenvalues as the radius increases with Algorithm \ref{alg:triangle} (first row) and Algorithm \ref{alg:meshless} (second row). The $x$-axis indicates the length of the lattice interval, and the $y$-axis indicates the computed eigenvalues. As the shell becomes thicker, the eigenvalues on the grid converge to eigenvalues on the smooth $S^2$ sooner. Fig.~A3 and Fig.~A4 display values of eigenvectors on $S^2$.

\setcounter{figure}{2}
\begin{figure}[!htb]
\includegraphics[width=0.3\textwidth]{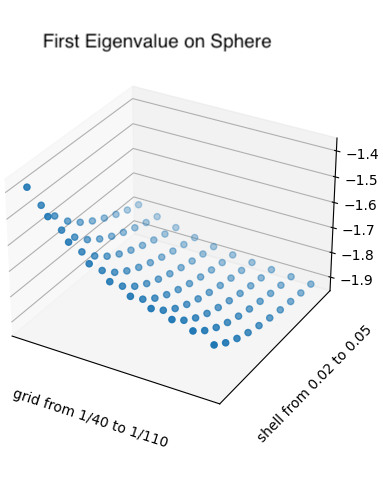}\hskip 3mm
\includegraphics[width=0.3\textwidth]{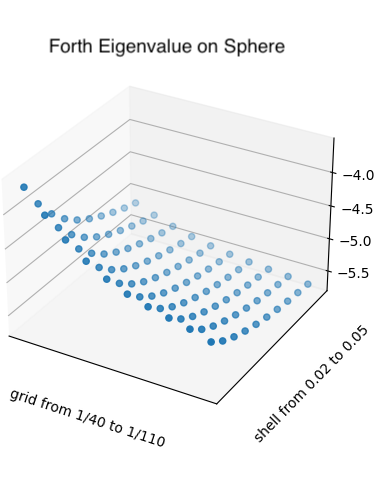}\hskip 3mm
\includegraphics[width=0.3\textwidth]{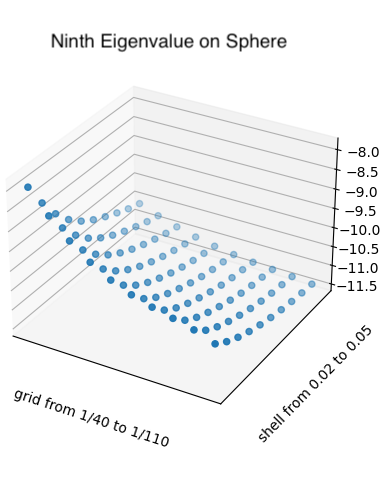}\hskip 3mm
\justify 
\caption{The first 3 non-zero eigenvalues (excluding multiplicity) on $S^2$ represented as a triangle mesh in 3D.
\textbf{Left:} The first eigenvalue $-\ell(\ell+1) = -2$ on $S^2$. \textbf{Middle:} The second eigenvalue $-\ell(\ell+1)=-6$ on $S^2$ with multiplicity $2\ell + 1$ for $\ell = 2$. \textbf{Right:} The third eigenvalue $-\ell(\ell+1)=-12$ on $S^2$ with multiplicity $2\ell + 1$ for $\ell = 3$.}\label{fig:ev_triangle_3D}
\end{figure}

\setcounter{figure}{3}
\begin{figure}[!htb]
\includegraphics[width=0.3\textwidth]{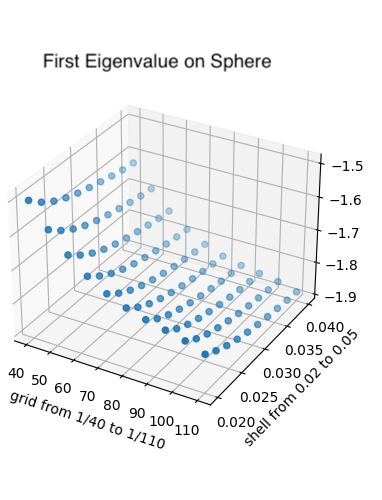}\hskip 5mm
\includegraphics[width=0.3\textwidth]{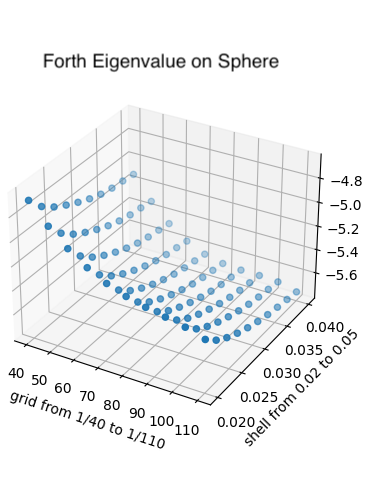}\hskip 5mm
\includegraphics[width=0.3\textwidth]{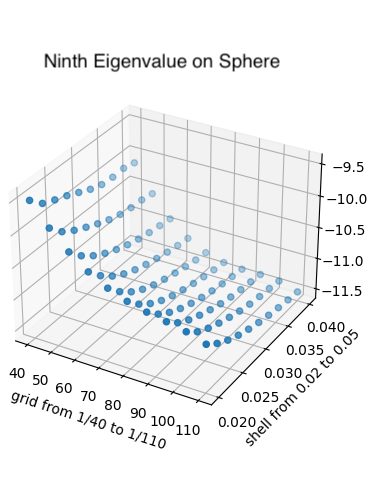}
\justify 
\caption{The first ten non-zero eigenvalues and multiplicity on $S^2$ represented as a point cloud in 3D. \textbf{Left:} The first eigenvalue $-\ell(\ell+1) = -2$ on $S^2$. \textbf{Middle:} The second eigenvalue $-\ell(\ell+1)=-6$ on $S^2$ with multiplicity $2\ell + 1$ for $\ell = 2$. \textbf{Right:} The third eigenvalue $-\ell(\ell+1)=-12$ on $S^2$ with multiplicity $2\ell + 1$ for $\ell = 3$.}\label{fig:ev_meshless_3D}
\end{figure}

\begin{figure}[!htb]
\hskip -8mm
\includegraphics[width=42mm]{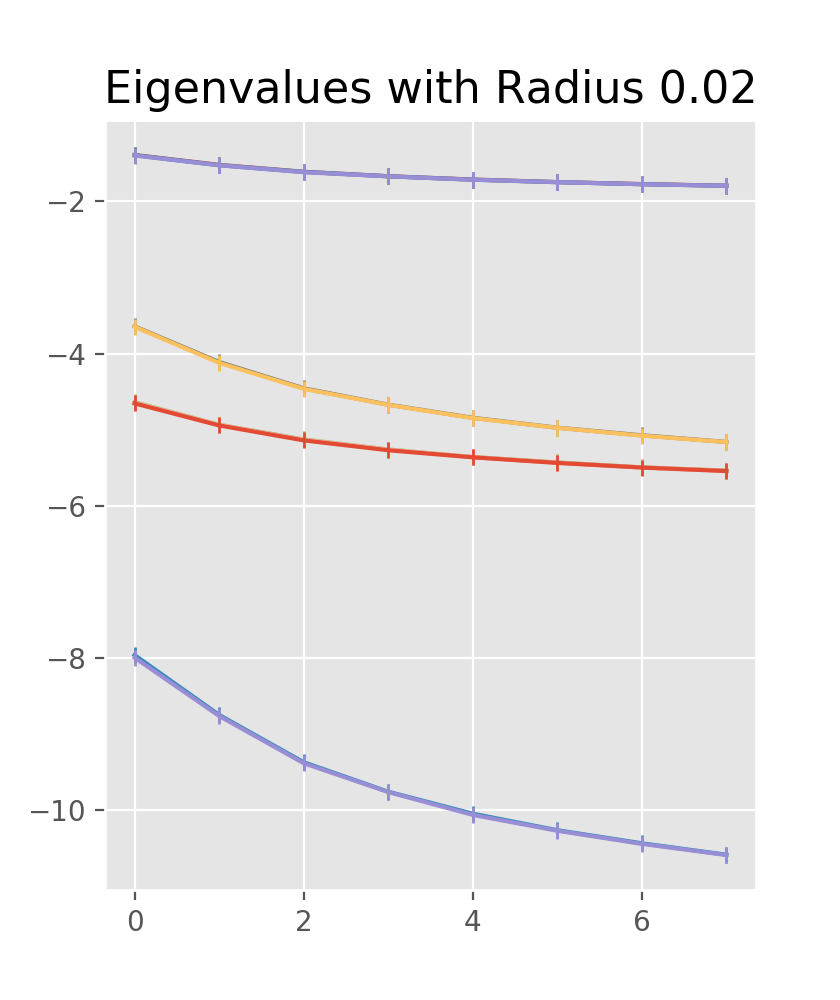}\hskip -3mm
\includegraphics[width=42mm]{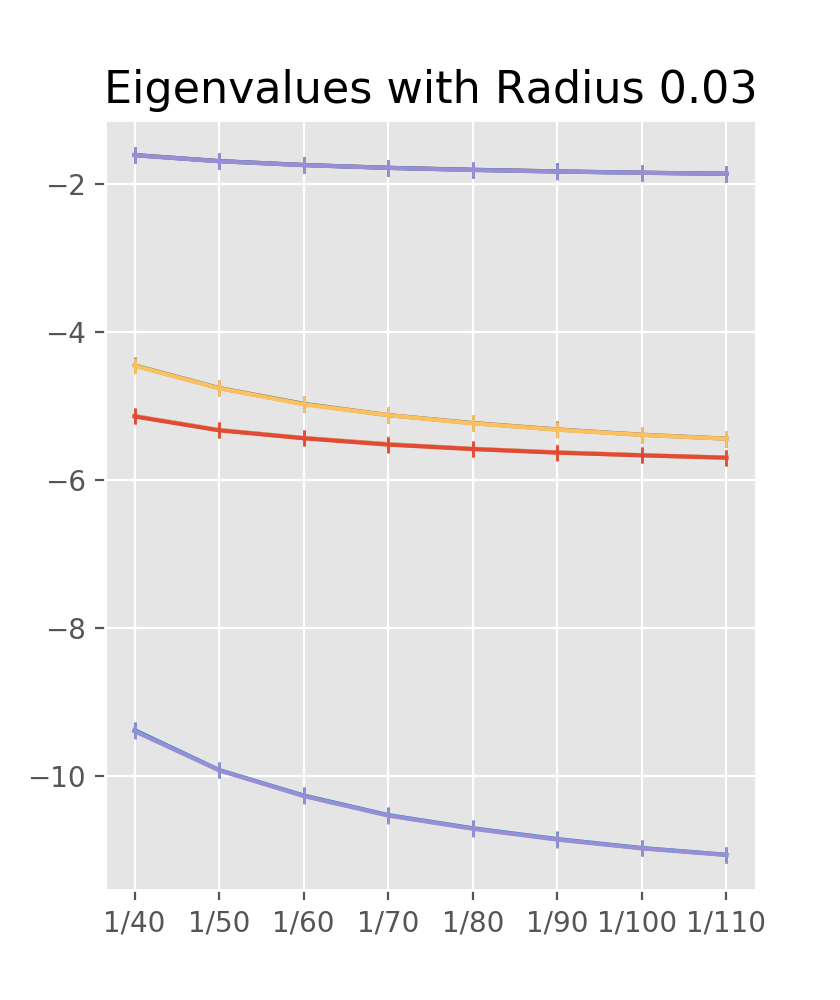}\hskip -3mm
\includegraphics[width=42mm]{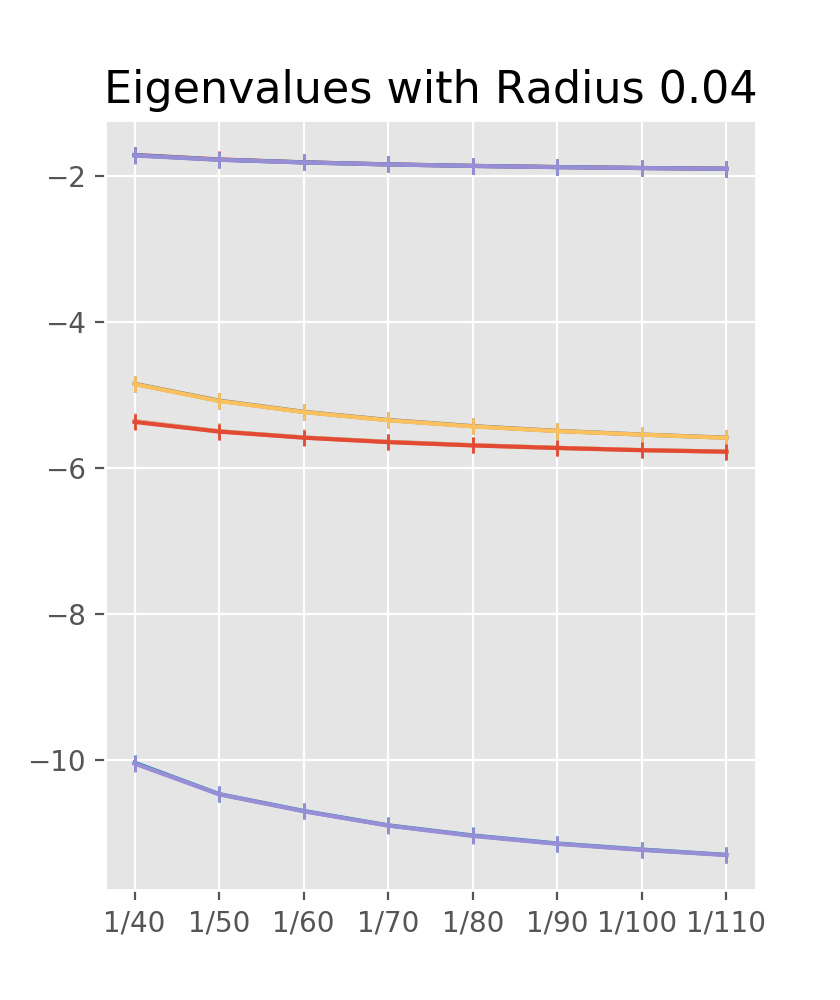}\hskip -3mm
\includegraphics[width=42mm]{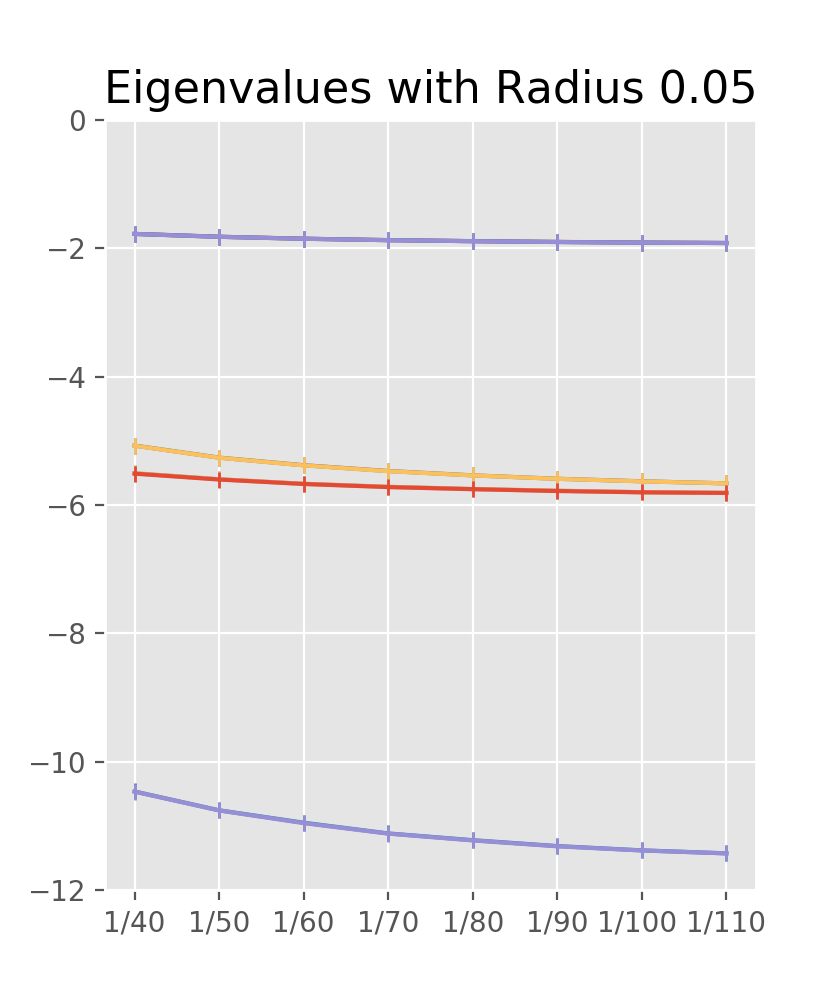}\\

\hskip -8mm
\includegraphics[width=42mm]{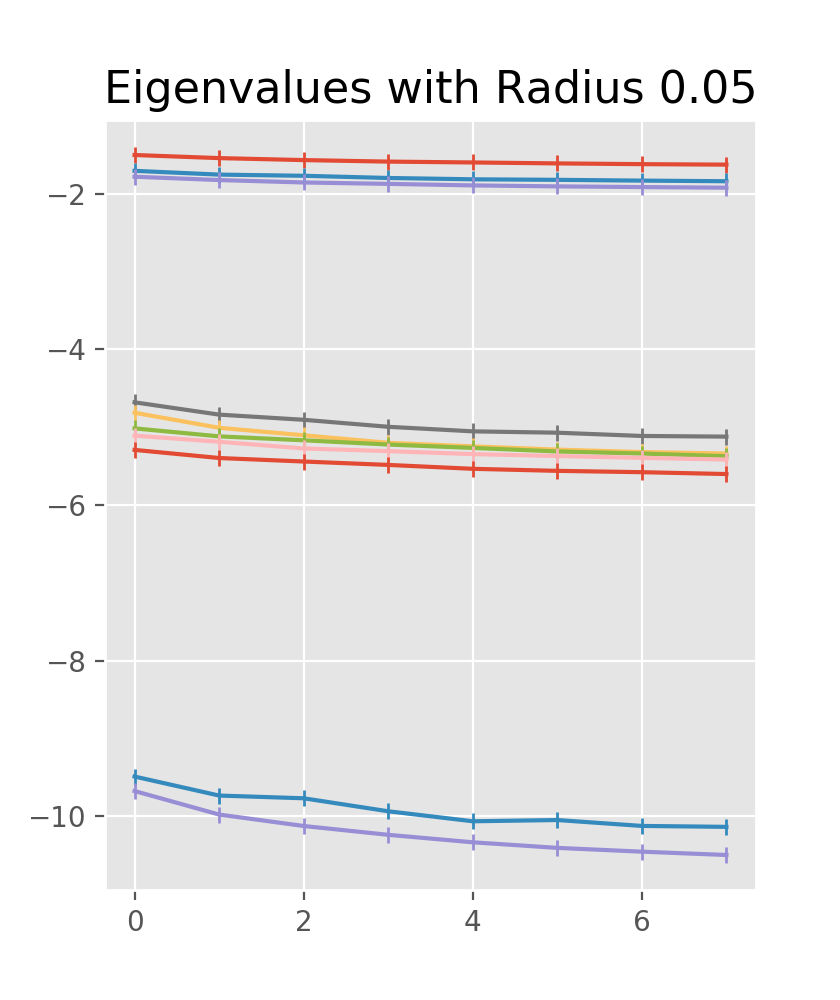}\hskip -3mm
\includegraphics[width=42mm]{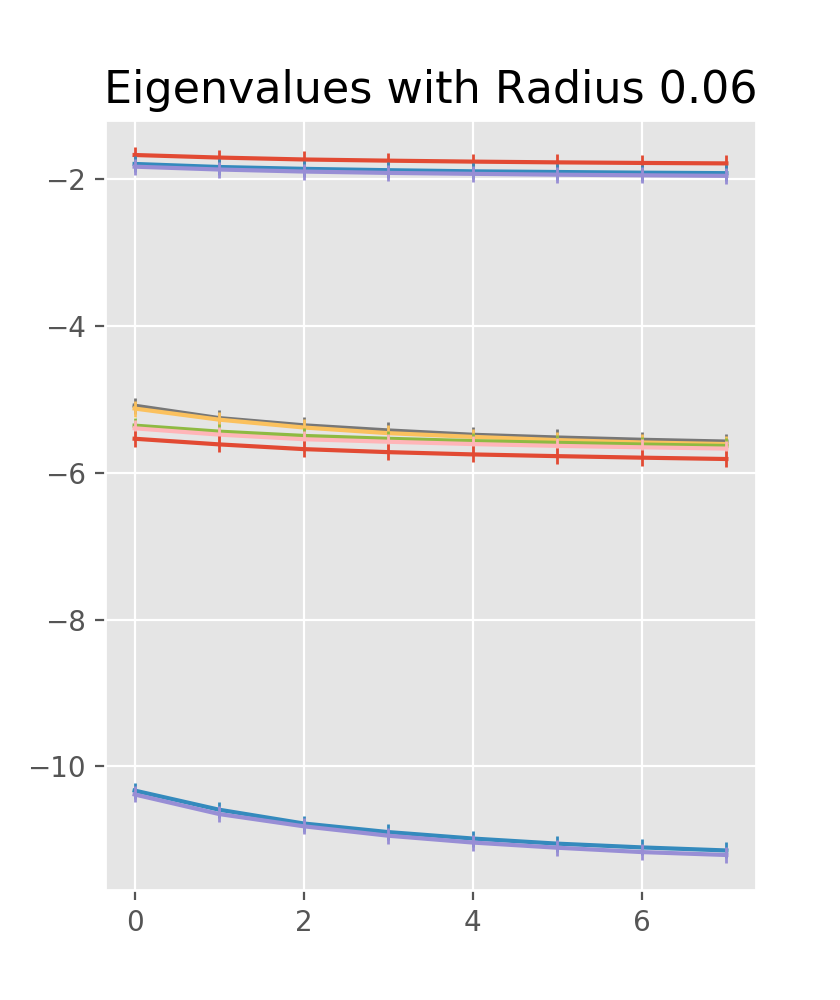}\hskip -3mm
\includegraphics[width=42mm]{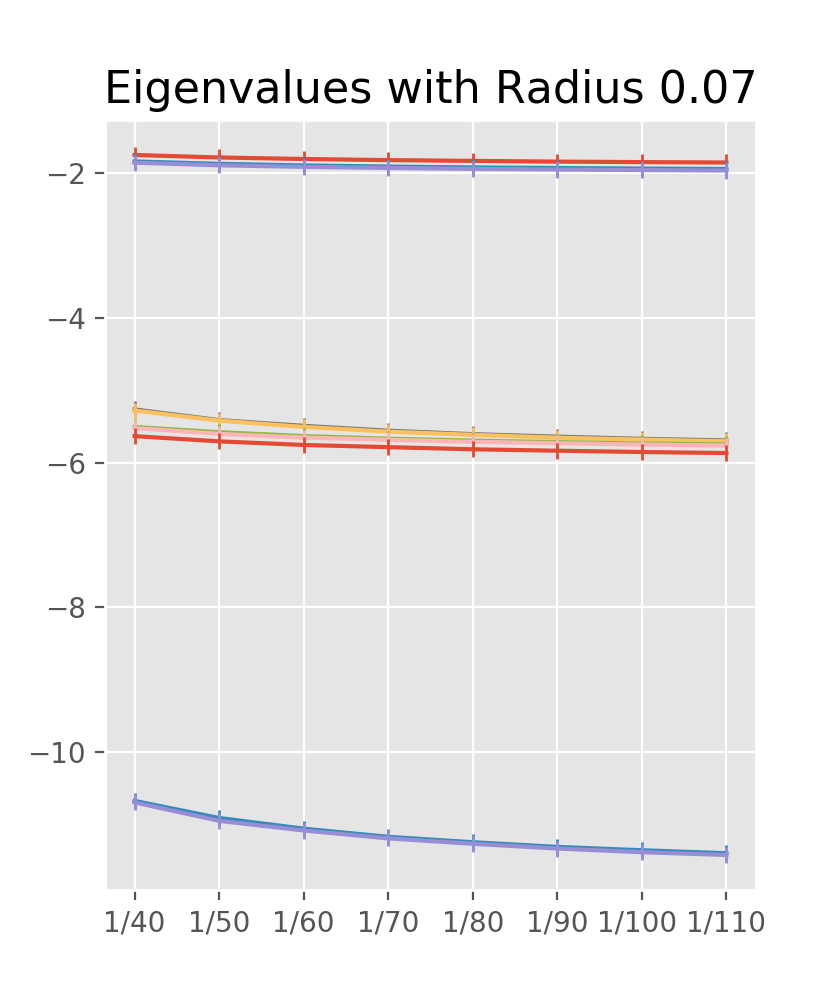}\hskip -3mm
\includegraphics[width=42mm]{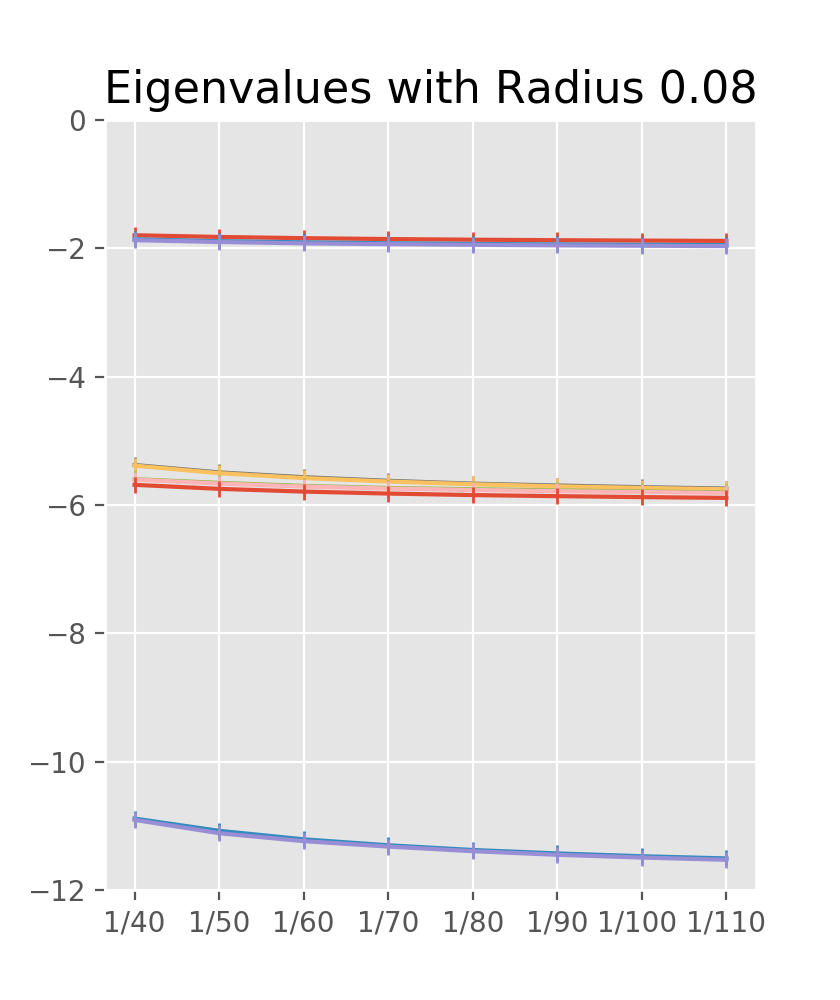}\hskip -3mm

\caption{Change of the first ten eigenvalues on $S^2$ represented as a point cloud as the shell radius increases with Algorithm \ref{alg:triangle} (first row) and Algorithm \ref{alg:meshless} (second row).}\label{fig:ev_meshless_radius}
\end{figure}

\subsection{Ellipsoids and Urakawa's Validation on Yau's Conjecture}

In this section, we report numerical experiments on the ten non-zero eigenvalues for five different ellipsoids using Algorithm \ref{alg:meshless}.
Table \ref{table:meshless_ellipsoid_1} reports the first ten non-zero surface eigenvalues on specified ellipsoids. In the table, CPU time measures the number of clock cycles summed across all threads that the Matlab process has used since Matlab started; the units are in seconds. 
All computations are performed on 170 GB memory allocation with 64 cores. 

In Jumonji and Urakawa's work on the eigenvalue problems of the Laplacian for bounded domains and embedded compact surfaces in the 3-dimensional Euclidean space \cite{jumonji2008eigenvalue}, they examined Yau's conjecture on the eigenvalues between embedded surfaces and their enclosed 3-dimensional bounded domains \cite{yau2000review}. In Jumonji and Urakawa's paper, they reported the first non-zero Laplacian eigenvalues of the ellipsoid with $a, b, c = 1,3,3$ are 0.3, 0.3, 0.4, 0.8, 0.8, 1.2, and 1.2, agreeing with our computation result, displayed in Table \ref{table:meshless_ellipsoid_1} in the fifth column.

\begin{figure}[!htb]
\begin{center}
\includegraphics[width=62mm]{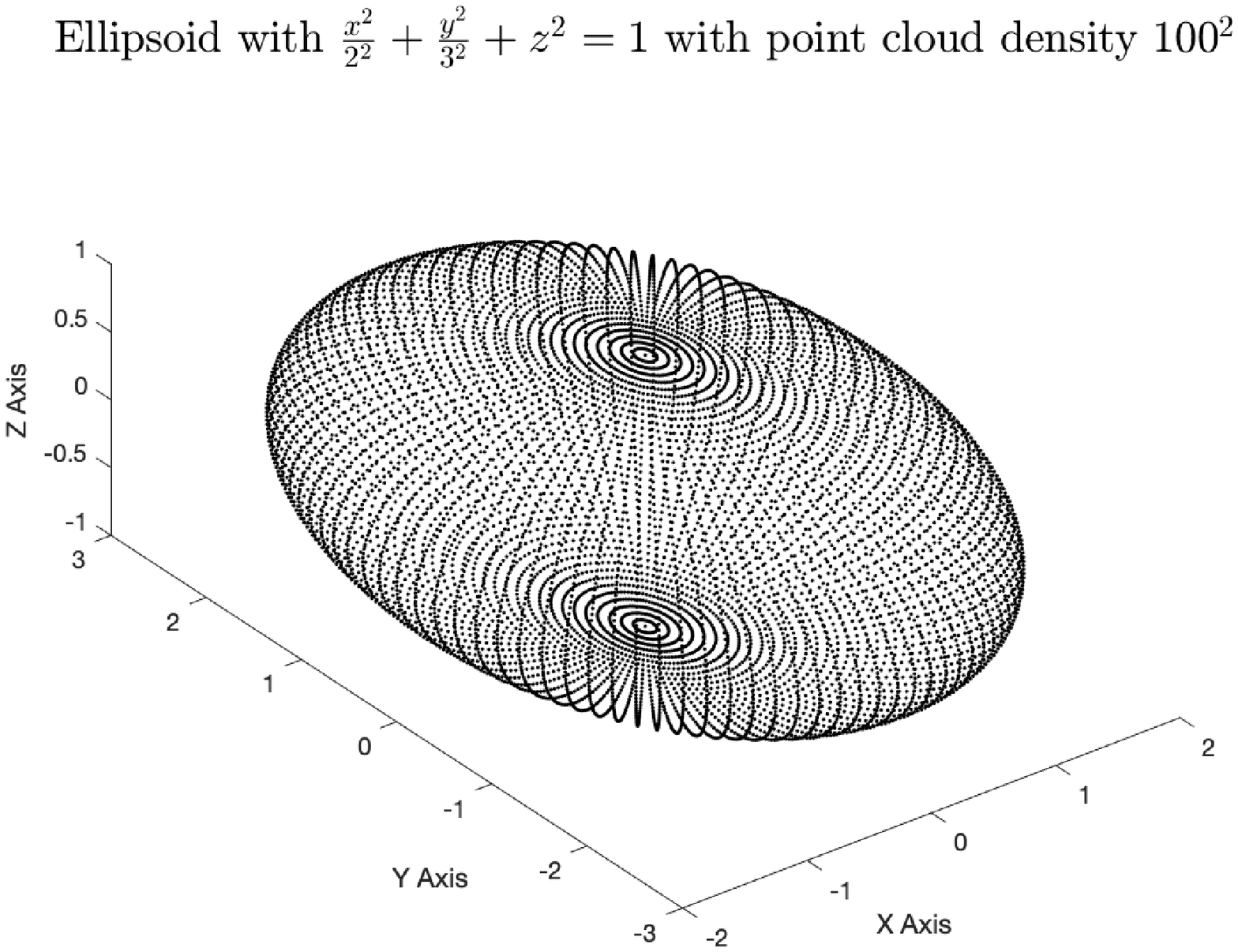}\hskip -3.6mm
\includegraphics[width=62mm]{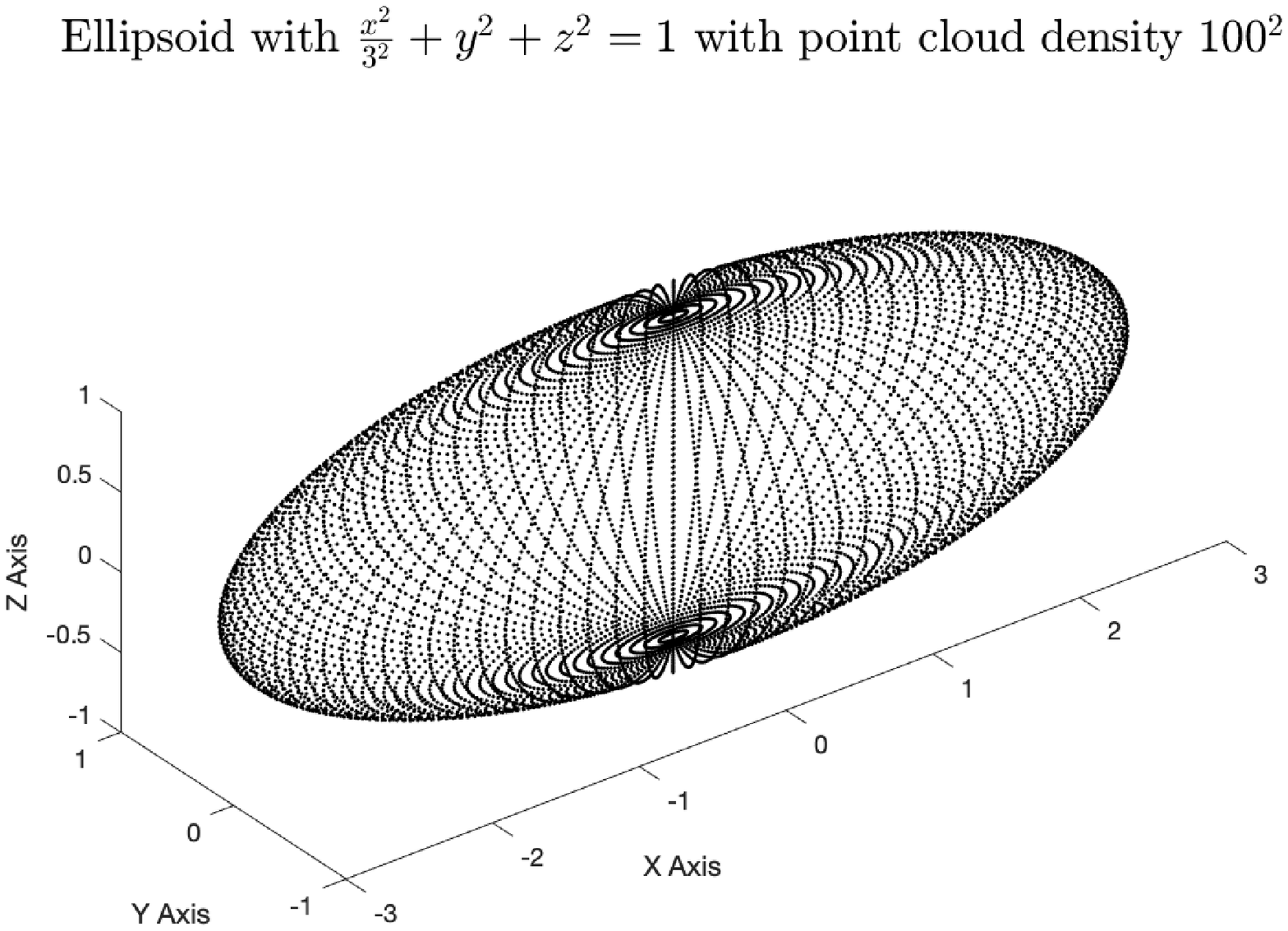}\hskip -3.6mm
\includegraphics[width=62mm]{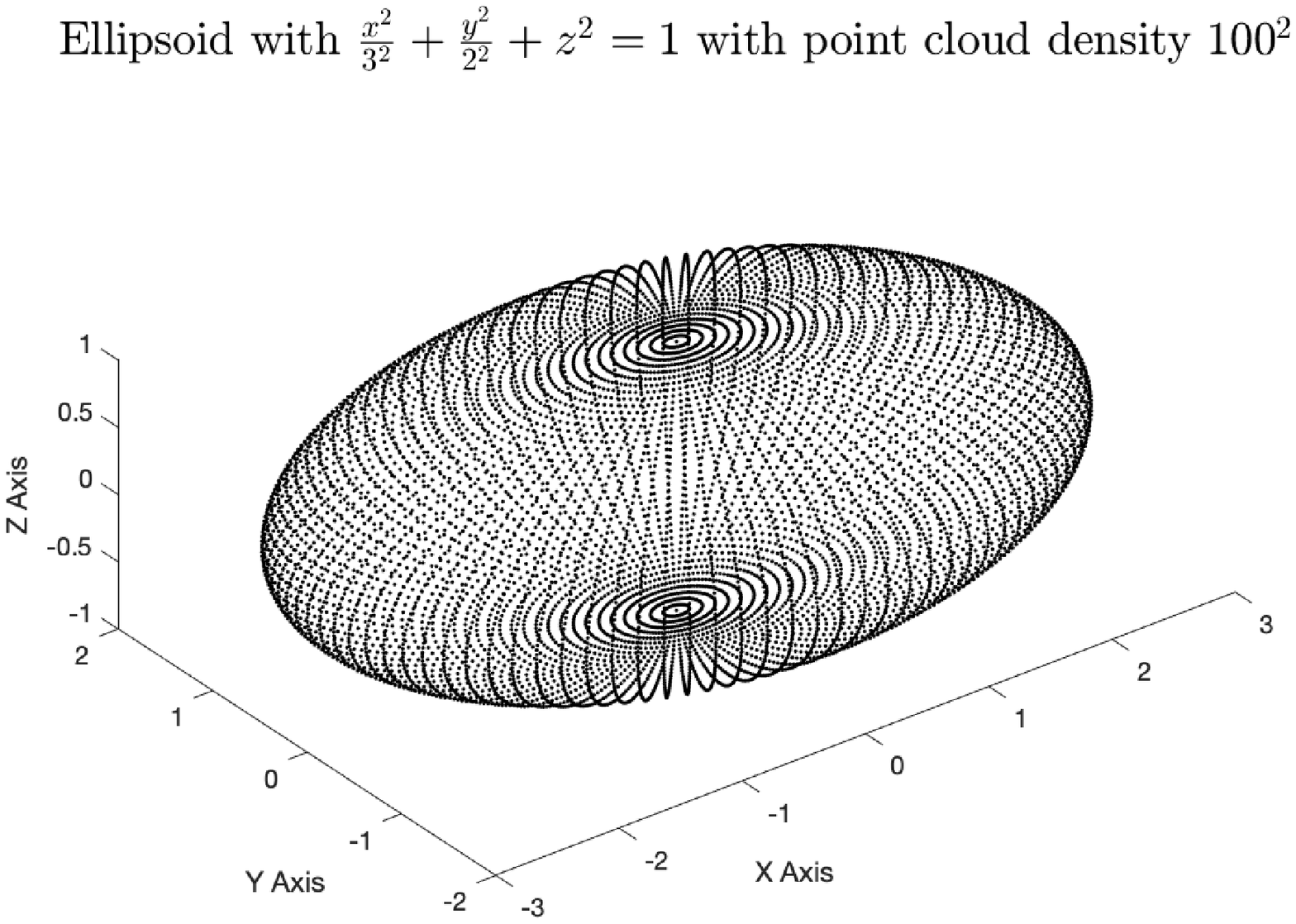}\hskip -3.6mm
\includegraphics[width=62mm]{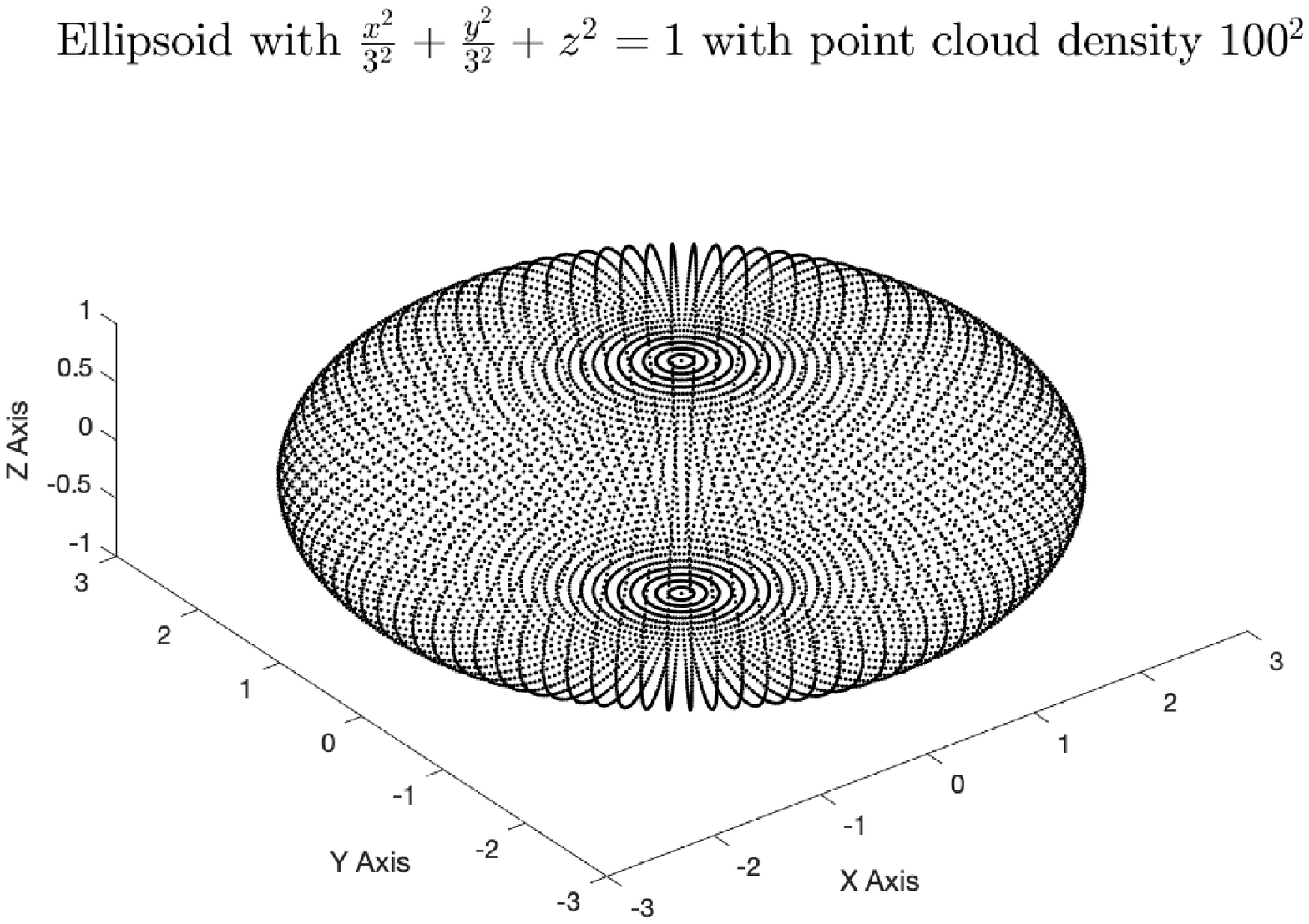}\hskip -3.6mm
\includegraphics[width=62mm]{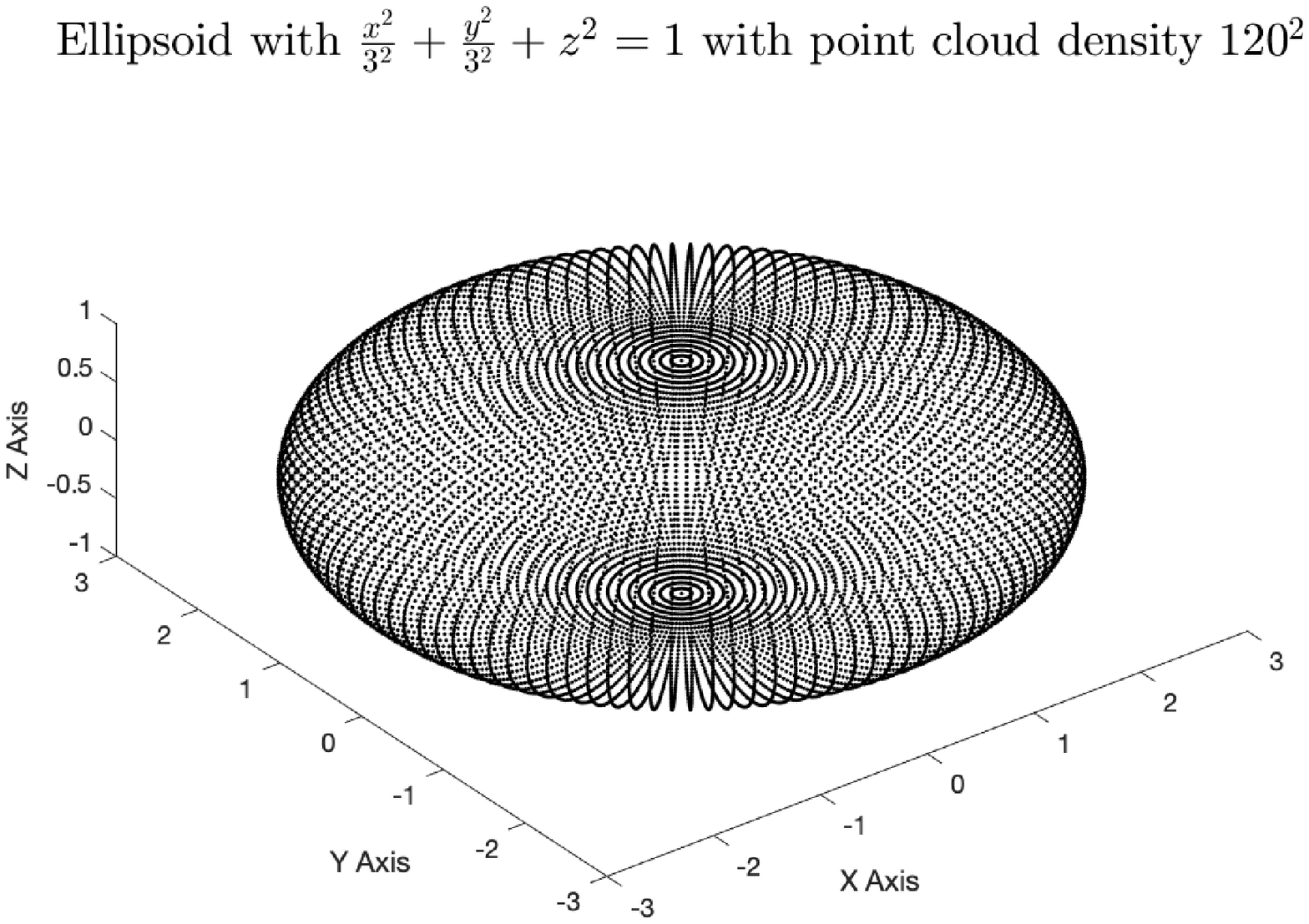}\hskip -3.6mm
\includegraphics[width=62mm]{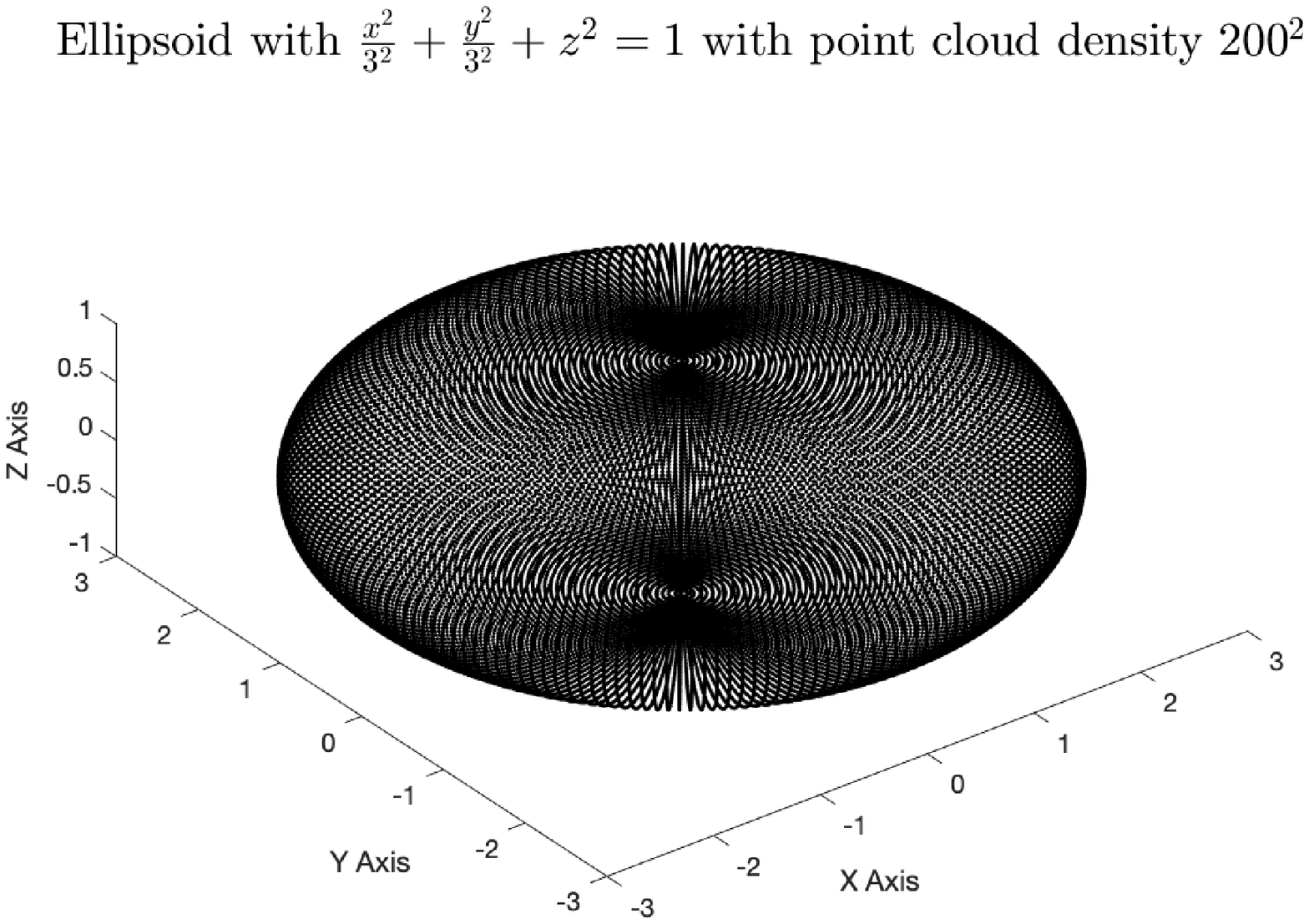}\hskip -3.6mm
\end{center}
\caption{Point cloud of ellipsoids.
\textbf{1:} The point cloud of the ellipsoid  $\displaystyle{x^2\over 2^2}+{y^2\over 3^2}+z^2=1$ with point density $100^2$.
\textbf{2:} The point cloud of the ellipsoid  $\displaystyle{x^2\over 3^2}+y^2+z^2=1$ with point density $100^2$.
\textbf{3:} The point cloud of the ellipsoid  $\displaystyle{x^2\over 3^2}+{y^2\over 2^2}+z^2=1$ with point density $100^2$.
\textbf{4:} The point cloud of the ellipsoid  $\displaystyle{x^2\over 3^2}+{y^2\over 3^2}+z^2=1$ with point density $100^2$.
\textbf{5:} The point cloud of the ellipsoid  $\displaystyle{x^2\over 3^2}+{y^2\over 3^2}+z^2=1$ with point density $120^2$.
\textbf{6:} The point cloud of the ellipsoid  $\displaystyle{x^2\over 3^2}+{y^2\over 3^2}+z^2=1$ with point density $200^2$.}\label{fig:ellipsoids}
\end{figure}

Fig.~\ref{fig:ellipsoids} exhibits the point cloud of ellipsoids $\displaystyle{x^2\over 2^2}+{y^2\over 3^2}+z^2=1$, $\displaystyle{x^2\over 3^2}+y^2+z^2=1$,  $\displaystyle{x^2\over 3^2}+{y^2\over 2^2}+z^2=1$, and $\displaystyle{x^2\over 3^2}+{y^2\over 3^2}+z^2=1$, respectively, with point density $100^2$, and the point cloud of the ellipsoid $\displaystyle{x^2\over 3^2}+{y^2\over 3^2}+z^2=1$ with point density $120^2$ and $200^2$. Here the point cloud density refers to the number of points in the point cloud that form the ellipsoid.
In particular, the eigenvalues of the ellipsoid $\displaystyle{x^2\over 2^2}+{y^2\over 3^2}+z^2=1$ should be the same as those of the ellipsoid $\displaystyle{x^2\over 3^2}+{y^2\over 2^2}+z^2=1$, since the latter is just an $SO(3)$ rotation of the first one, which can be observed accordingly from Table \ref{table:meshless_ellipsoid_1}.

\begin{table}[!htb]
\centering
\caption{The first ten non-zero surface eigenvalues by the meshless algorithm with $r = 0.11$ and $n=50$ on 64 Cores and 170GB Memory HPC.}\label{table:meshless_ellipsoid_1}
\begin{tabular}{|l|r|r|r|r|r|r|r|r|}
\hline
\hline
&\small{$\displaystyle{x\over 2^2}+{y\over 3^2}+z^2=1$}&\small{$\displaystyle{x^2\over 3^2}+y^2+z^2=1$}&\small{$\displaystyle{x^2\over 3^2}+{y^2\over 2^2} +z^2=1$}&\small{$\displaystyle{x^2\over 3^2}+{y^2\over 3^2} +z^2=1$} & \small{$\displaystyle{x^2\over 3^2}+{y^2\over 2} +z^2=1$}\\
\hline
\multirow{10}*{Eigenvalues}
&-0.315241&-0.253582&-0.315241&-0.308855&-0.3228\\
&-0.596258&-1.126602&-0.596258&-0.308855& -0.9139\\
&-0.662100&-1.230701&-0.662100&-0.457628& -0.9425\\
&-1.017389&-1.363820&-1.017389&-0.833177& -1.0820\\
&-1.243444&-2.203574&-1.243444&-0.833941& -1.7189\\
&-1.486521&-2.276262&-1.486521&-1.219271& -1.8059\\
&-2.045236&-2.336757&-2.045236&-1.219271& -2.2293\\
&-2.048874&-3.549628&-2.048874&-1.269028& -2.8636\\
&-2.060981&-3.584377&-2.060981&-1.562024& -3.0742\\
&-2.216862&-3.975534&-2.216862&-1.562024 & -3.1174\\
\hline
Elapsed Time &  736'' &  425'' & 795''& 1,341''&  639''\\
CPU Time & 5,043'' & 2,413'' & 5,440''&7,770'' & 4,057''\\
\hline
\hline
\end{tabular}
\label{tab1}
\end{table}


Table \ref{table:meshless_ellipsoid_4_equal_dist} compares the surface eigenvalues of the ellipsoid $\displaystyle{x^2\over 3^2}+{y^2\over 3^2} +z^2=1$ with various point cloud densities. Comparing to results reported in Jumonji and Urakawa, both the point cloud density and meshless algorithm parameters contribute to the convergence of eigenvalues. More specifically, greater point cloud density, thicker grid (greater shell thickness $r$), and denser grid (greater $n$) all lead to a more precise eigenvalue.

\begin{table}[!htb]
\centering
\caption{The first ten non-zero surface eigenvalues by the meshless algorithm with 64 Cores and 170GB Memory with point cloud density $100^2$, $120^2$, $200^2$, $500^2$, and $1000^2$.}\label{table:meshless_ellipsoid_4_equal_dist}
\begin{tabular}{|l|r|r|r|r|r|r|r|r|}
\hline
\hline\multicolumn{9}{|c|}{$\displaystyle{x^2\over 3^2}+{y^2\over 3^2} +z^2=1$}\\
\hline
&\multicolumn{3}{c|}{Point Cloud Density $100^2$}&\multicolumn{2}{c|}{Point Cloud Density $120^2$}&\multicolumn{3}{c|}{Point Cloud Density $200^2$}\\
\hline
&$r = 0.05$ & $r = 0.05$ & $r = 0.05$ &  $r = 0.05$ &  $r = 0.05$  &$r = 0.05$ & \multicolumn{2}{c|}{$r = 0.05$ } \\
&$n = 50$ & $n = 75$& $n = 100$&$n = 50$ & $n = 75$&$n = 50$ &  \multicolumn{2}{c|}{$n = 75$}\\
\hline
\multirow{10}*{Eigenvalues}
&-0.173892&-0.178447&-0.180789&-0.299862&-0.303826&-0.299862& \multicolumn{2}{c|}{-0.303826} \\
&-0.173892&-0.178568&-0.180789&-0.299862&-0.303885&-0.299862 & \multicolumn{2}{c|}{-0.303885}\\
&-0.263055&-0.273016&-0.276887&-0.415050 &-0.434217&-0.415050& \multicolumn{2}{c|}{-0.434217}\\
&-0.265516&-0.273873&-0.277627&-0.793171&-0.806634 &-0.793171& \multicolumn{2}{c|}{-0.806634}\\
&-0.284125&-0.301827&-0.312457&-0.795942&-0.807992 &-0.795942& \multicolumn{2}{c|}{-0.807992}\\
&-0.318465&-0.331424&-0.336788&-1.145733&-1.177993 &-1.145733& \multicolumn{2}{c|}{-1.177993}\\
&-0.318465&-0.332699&-0.336788&-1.145733&-1.178085 &-1.145733& \multicolumn{2}{c|}{-1.178085}\\
&-0.347359&-0.368679&-0.374133 &-1.193512&-1.238377&-1.193512& \multicolumn{2}{c|}{-1.238377}\\
&-0.365683&-0.376422&-0.381585&-1.471670&-1.496612&-1.471670& \multicolumn{2}{c|}{-1.496612}\\
&-0.383497&-0.399842&-0.407944&-1.471670&-1.497865&-1.471670 & \multicolumn{2}{c|}{-1.497865}\\
\hline
Elapsed Time & 63'' &  305'' & 1,262'' & 173'' & 1,137'' & 173''& \multicolumn{2}{c|}{1,116''}\\
CPU Time & 249'' & 1,183'' & 4,210'' & 949'' & 7,255'' & 938''& \multicolumn{2}{c|}{7,538''}\\
\hline\hline
&\multicolumn{3}{c|}{Point Cloud Density $500^2$}&\multicolumn{3}{c|}{Point Cloud Density $1000^2$}& $500^2$&  $1000^2$\\
\hline
&$r = 0.05$ & $r = 0.05$ & $r = 0.05$ &  $r = 0.05$ &  $r = 0.05$  &$r = 0.05$ &   $r = 0.08$&  $r = 0.08$   \\
&$n = 50$ & $n = 75$& $n = 100$&$n = 50$ & $n = 75$&$n = 100$ & $n = 50$ & $n = 50$\\
\hline
\multirow{10}*{Eigenvalues} &
-0.303699&-0.308881&-0.311363&-0.302948&-0.308269&-0.310890&-0.309011&-0.308697\\
&-0.303699&-0.308910&-0.311363&-0.302948&-0.308324&-0.310890&-0.309011&-0.308697\\
&-0.447712&-0.466251&-0.475603&-0.451377&-0.469358&-0.478537&-0.471729&-0.472977\\
&-0.822690&-0.844795&-0.855671&-0.821883&-0.844126&-0.855192&-0.846485&-0.845983\\
&-0.846778&-0.860377&-0.866436&-0.847143&-0.861282&-0.867818&-0.862488&-0.862434\\
&-1.201132&-1.233425&-1.249930&-1.205294&-1.237182&-1.253610&-1.242967&-1.244217\\
&-1.201132&-1.233541&-1.249930&-1.205294&-1.237580&-1.253610&-1.242967&-1.244217\\
&-1.204181&-1.245586&-1.265774&-1.206545&-1.246583&-1.266857&-1.250757&-1.251632\\
&-1.582467&-1.619629&-1.637757&-1.583844&-1.621378&-1.640358&-1.625266&-1.625361\\
&-1.582467&-1.620202&-1.637757&-1.583844&-1.622490&-1.640358&-1.625266&-1.625361\\
\hline
Elapsed Time & 258'' & 1,516'' & 7,894''&290''& 1,885'' & 8,239'' & 705'' & 898''\\
CPU Time & 1,346'' & 9,588'' & 59,145''& 1,331'' & 11,290'' & 59,262'' & 4,390'' & 4,760''\\
\hline
\hline
\end{tabular}
\label{tab1}
\end{table}

\subsubsection{Tetrahedron}
We present the Laplacian eigenvalues of a 3-sided tetrahedron and a 4-sided tetrahedron, shown in Fig.~A1. The results are displayed in Table \ref{tab:3-4-tetra}, where \it Surface \rm column indicates eigenvalues computed by the classical algorithm for triangle mesh by Reuter \cite{reuter2009discrete}, the \it Mesh Vertex \rm column represents the eigenvalues computed by meshless algorithm using the point cloud consists of the vertex set from a triangle mesh, and the \it Point Cloud \rm column represent eigenvalues computed by meshless algorithm (Algorithm 3) using a generated point cloud of respective surface. A discussion about the analytic solution to a tetrahedron could be found in \cite{greif2018spectrum}.

\begin{table}[!htb]
\centering
\centering\caption{The first ten non-zero surface eigenvalues of the 3-sided regular tetrahedron and 4-sided regular tetrahedron by triangle surface algorithm and the meshless algorithm with point cloud density $100^2$ and shell thickness 0.05.}\label{tab:tetra}
\begin{tabular}{|l|c|r|r|r|c|r|r|r|}
\hline
&\multicolumn{3}{c|}{3-sided Regular Tetrahedron}&\multicolumn{3}{c|}{4-sided Regular Tetrahedron}\\
\hline
&\multirow{2}*{Surface}&\multicolumn{2}{c|}{Meshless algorithm}&\multirow{2}*{Surface}&\multicolumn{2}{c|}{Meshless algorithm}\\
\cline{3-4}\cline{6-7}
& &Point Cloud & Mesh Vertex
& & Original & Mesh Vertex
\\
\hline
\multirow{10}*{Eigenvalues}
&-3.6214&-3.3574& -3.3450  &-4.3830&-4.0153  &-4.2557\\
&-3.6214&-3.4153& -3.3782  &-4.3830&-4.2599   &-4.3071\\
&-5.8424&-4.8823& -5.2485  &-4.3830&-4.3330   &-4.3378\\
&-11.1263&-10.5358& -10.4637  &-13.1278&-12.8418 &-13.2930\\
&-11.1263&-10.6727& -10.5752  &-13.1278&-13.1429 &-13.3975\\
&-17.4897&-15.0514& -15.6586  &-13.1278&-13.2677  &-13.5373\\
&-18.1964&-16.3376& -16.9002  &-17.4897&-16.2846  &-16.5812\\
&-18.1964&-16.5558& -17.0233  &-17.4897&-17.0448 &-17.4354\\
&-23.2946&-20.1906& -20.9053  &-17.4897&-17.1762  &-17.5524\\
&-23.7987&-23.5115&  -23.2507  &-17.4897&-30.5332  &-30.2258\\
\hline
Elapsed Time &0.84'' & 157'' & 118''  &0.72''& 249''&217''\\
CPU Time &7.88'' & 878'' & 624''   &14''&  1,445''&1,071''\\
\hline
\hline
\end{tabular}
\label{tab:3-4-tetra}
\end{table}

\subsubsection{Cube}


We present eigenvalues of a cube with vertices at $(\pm 1, \pm 1, \pm 1)$. The computed eigenvalues are presented in Table \ref{tab:cube}. 
 The triangle mesh used for the cube is displayed in Fig.~A2  
on the left, and we also used a separate point cloud surface of the cube.
 In Table \ref{tab:cube}, the forth column, \it Mesh Vertex, \rm represents the eigenvalues computed by meshless algorithm (Algorithm 3) using the point cloud consists of the vertex set from a triangle mesh. The first three columns are computed eigenvalues of a cube represented by a generated point cloud.
\begin{table}[!htb]
\centering
\centering\caption{The first ten non-zero surface eigenvalues of a cube with vertices at $(\pm 1, \pm 1, \pm 1)$ by the meshless algorithm with point cloud density $100^2$ and shell thickness 0.05.}\label{tab:cube}
\begin{tabular}{|l|c|c|c|c|}
\hline
&30 points on each edge&50 points on each edge&Mesh Vertex&Surface\\
\hline
\multirow{10}*{Eigenvalues}
&   -0.9416 &-1.0571&-0.9781&-1.0412\\
&   -0.9416&-1.0571& -0.9781&-1.0412\\
&   -0.9417&-1.0571&-0.9781&-1.0412\\
&   -2.5737&-2.8859&-2.6651&-2.8728\\
&   -2.5737&-2.8859&-2.6651&-2.8728\\
&   -2.5737&-2.8859&-2.6651&-2.8728\\
&   -3.2527&-3.6592&-3.3990&-3.5307\\
&   -3.2528&-3.6592&-3.3990&-3.5307\\
&   -4.4080&-4.9393&-4.5610&-4.9308\\
&   -5.8522&-6.5749&-6.0938&-6.4163\\
\hline
Elapsed Time &1,066''&1,768''&1,492''&0.79''\\
CPU Time &6,320'' &10,987''&8,036''&740''\\
\hline
\hline
\end{tabular}
\end{table}

\subsubsection{Cone}
In this section, we report the computed eigenvalues using the classical surface algorithm and meshless algorithm (Algorithm 3); the results are reported in Table \ref{tab:cone_density}, where the first three columns report eigenvalues computed by the meshless algorithm, and the last column reports eigenvalues computed by the surface algorithm. In particular, the first two columns report eigenvalues computed by meshless algorithm using generated point cloud, and the third column reports eigenvalues computed by meshless algorithm using the point cloud extracted from the triangle mesh, as exhibited in Fig.~A2 
on the right.
We also display values of eigenvectors on this cone in Fig.~A5 and Fig.~A6.

\begin{table}[!htb]
\centering
\centering\caption{The first ten non-zero surface eigenvalues of a cone with vertex at the origin and base at $z=1$ with radius = 1 by the meshless algorithm with point cloud density $100^2$ and shell thickness 0.05.}\label{tab:cone_density}
\begin{tabular}{|l|c|c|c|c|}
\hline
&20 points on each edge&50 points on each edge&Mesh Vertex&Surface\\
\hline
\multirow{10}*{Eigenvalues}
&-1.3725 &-2.3616&-3.0001&-3.0376\\
&   -2.8905&-3.2116&-3.0001&-3.0377\\
&   -3.0874&-3.9597&-4.1453&-3.8591\\
&   -5.6597&-6.4916&-8.4367&-8.5843\\
&   -5.7116&-8.6735&-8.5087&-8.5854\\
&   -6.5111&-10.7776&-10.6293&-10.3370\\
&   -7.2630&-11.1218&-12.3705&-11.4486\\
&   -7.8935&-11.9478&-12.3705&-11.4491\\
&  -10.0452&-13.4958&-16.2407&-16.4751\\
&  -10.3218&-14.0009&-16.2407&-16.4761\\
\hline
Elapsed Time &147''&324''&380''&2.02''\\
CPU Time &705'' & 1,815''&2,006''&2,049''\\
\hline
\hline
\end{tabular}
\end{table}

\subsection{Computed First Ten Non-zero Eigenvalues of Genus 1 Surface}

Table \ref{table:meshless_genus1} reports the first ten non-zero surface eigenvalues by the meshless algorithm of a genus 1 surface.\\

\begin{minipage}{.25\textwidth}
\includegraphics[width=0.7\textwidth]{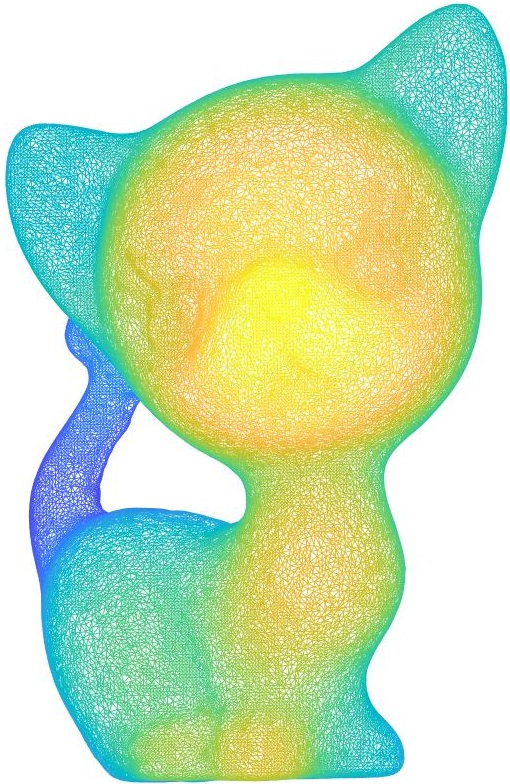}
    \end{minipage}%
    \begin{minipage}{0.75\textwidth}
\captionof{table}{The first ten non-zero surface eigenvalues by the meshless algorithm with 64 Cores and 170GB.}\label{table:meshless_genus1}
\begin{tabular}{|l|r|r|r|r|r|}
\hline
&$r = 0.05$ & $r = 0.05$ & $r = 0.08$ &  $r = 0.08$ &  $r = 0.11$ \\
&$n = 75$ & $n = 100$&$n = 75$ & $n = 75$&$n = 100$ \\
\hline
\multirow{10}*{Eigenvalues}
&-0.617180&-0.628299&-0.646848&-0.653546&-0.672847\\
&-2.079415&-2.112728&-2.152381&-2.173507&-2.200838\\
&-2.494772&-2.533026&-2.582687&-2.606614&-2.640699\\
&-3.088792&-3.136877&-3.198115&-3.227173&-3.273606\\
&-3.391403&-3.436752&-3.511123&-3.544625&-3.605186\\
&-4.661078&-4.721061&-4.821988&-4.861155&-4.933898\\
&-5.473066&-5.574882&-5.698128&-5.751677&-5.834958\\
&-6.534291&-6.644320&-6.756330&-6.824672&-6.883822\\
&-6.632011&-6.733570&-6.904994&-6.969889&-7.129169\\
&-7.239121&-7.338298&-7.541480&-7.603741&-7.796873\\
\hline
Elapse Time & 216''&882''&686''&3,241''&1,421''
\\
CPU Time & 1,166''&5,743''&4,543''&23,234''&9,871''
\\
\hline\hline
\end{tabular}
    \end{minipage}

\subsection{Computed First Ten Non-zero Eigenvalues of Genus 2 Surface}

Table \ref{table:meshless_genus2} reports the first ten non-zero surface eigenvalues by the meshless algorithm of a genus 2 surface.\\

\begin{minipage}{.4\textwidth}
\includegraphics[width=0.5\textwidth]{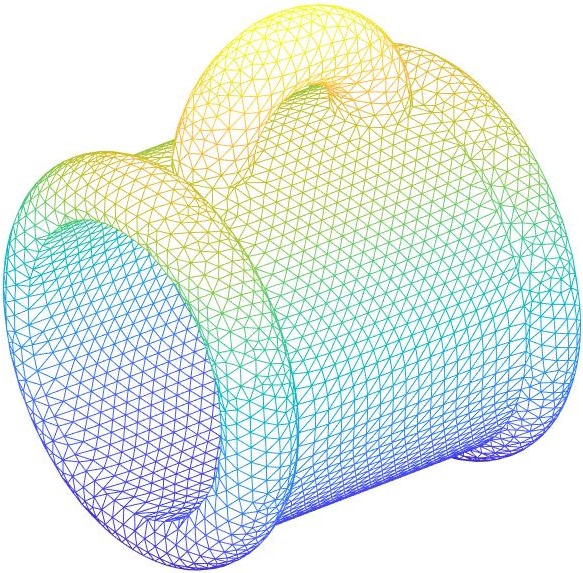}
    \end{minipage}%
    \begin{minipage}{0.5\textwidth}
%
%
\captionof{table}{The first ten non-zero surface eigenvalues by the meshless algorithm with 64 Cores and 170GB.}\label{table:meshless_genus2}
\begin{tabular}{|l|r|r|r|r|r||r|r|r|r|r|r|}
\hline
&$r = 0.05$ & $r = 0.05$ & $r = 0.08$ \\
&$n = 75$& $n = 100$&  $n = 75$\\
\hline
\multirow{10}*{Eigenvalues}
&-0.750970&-0.764963&-0.907334\\
&-0.955846&-0.973025&-1.124415\\
&-1.167671&-1.182242&-1.628818\\
&-1.419964&-1.429386&-1.667057\\
&-1.990223&-2.017348&-2.518728\\
&-2.292042&-2.315880&-2.733560\\
&-2.352482&-2.393655&-2.819258\\
&-2.386239&-2.414864&-3.015961\\
&-2.786546&-2.824579&-3.486652\\
&-3.308398&-3.375575&-4.203104\\
\hline
Elapse Time &199''&626''&1,449''
\\
CPU Time & 1,041''&3,946''&11,563''
\\
\hline\hline
\end{tabular}
\end{minipage}

\subsection{Face}
\begin{figure}[!htb]
\includegraphics[width=55mm]{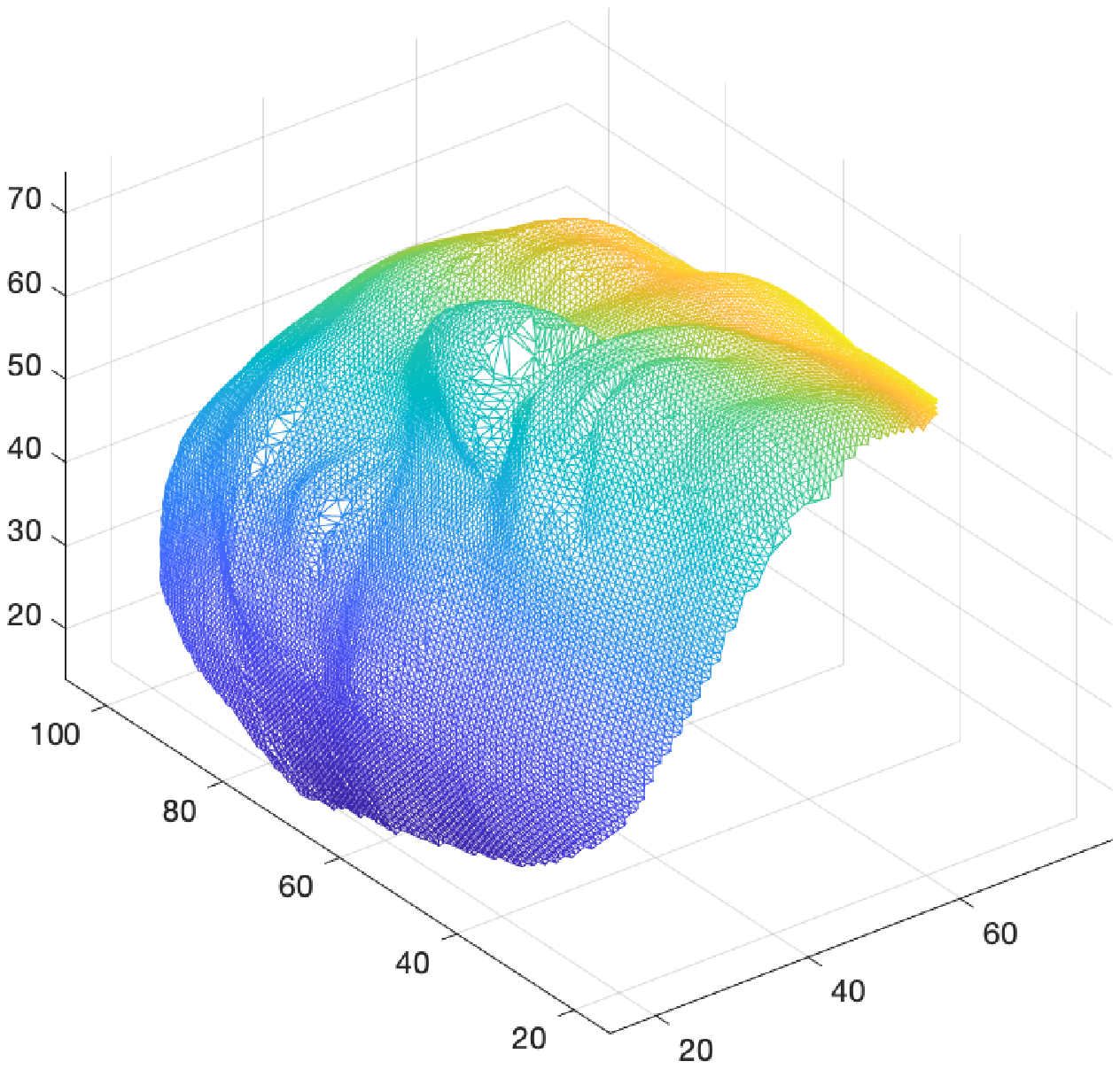}\hskip -12mm
\includegraphics[width=55mm]{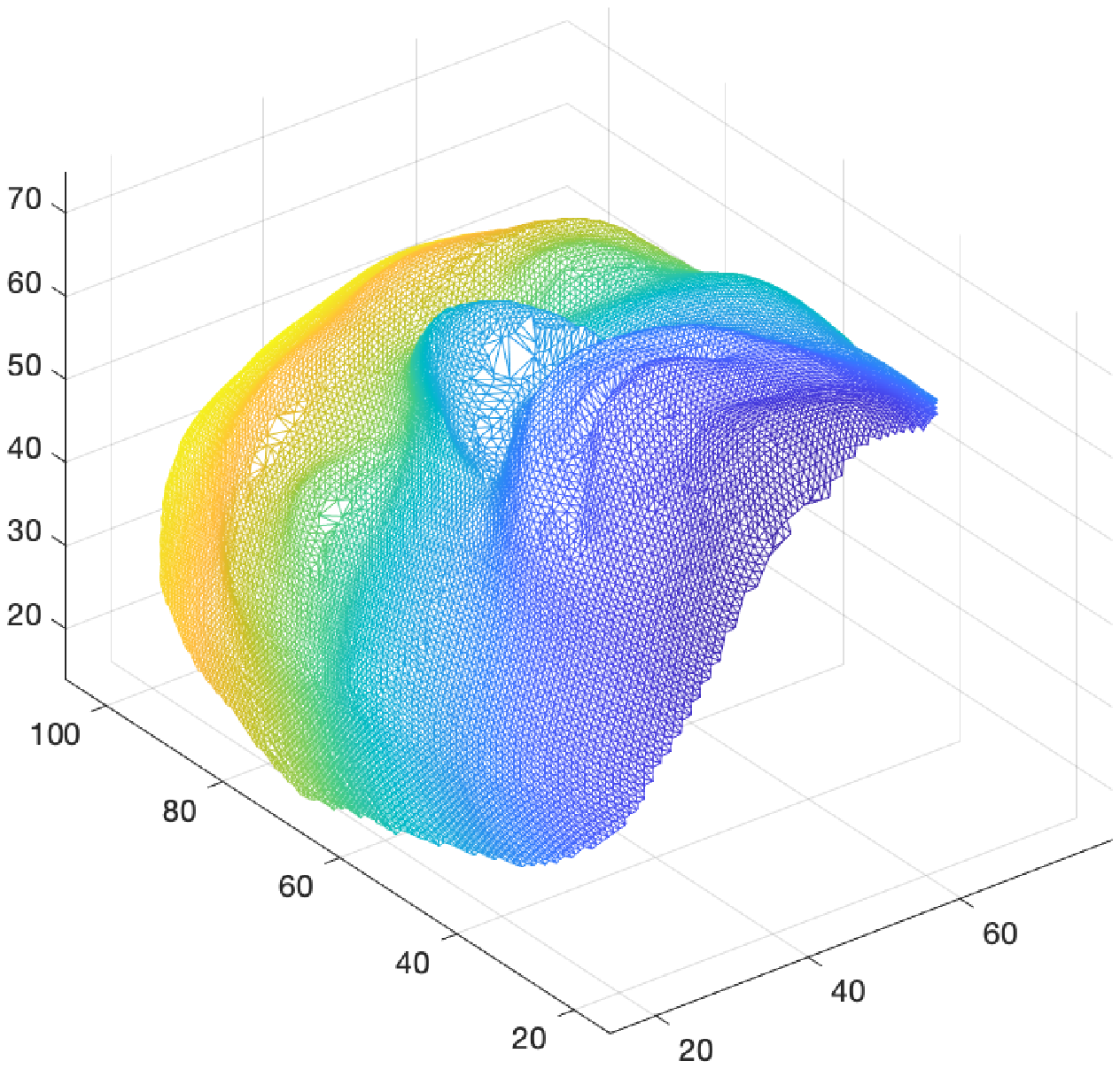}\hskip -12mm
\includegraphics[width=55mm]{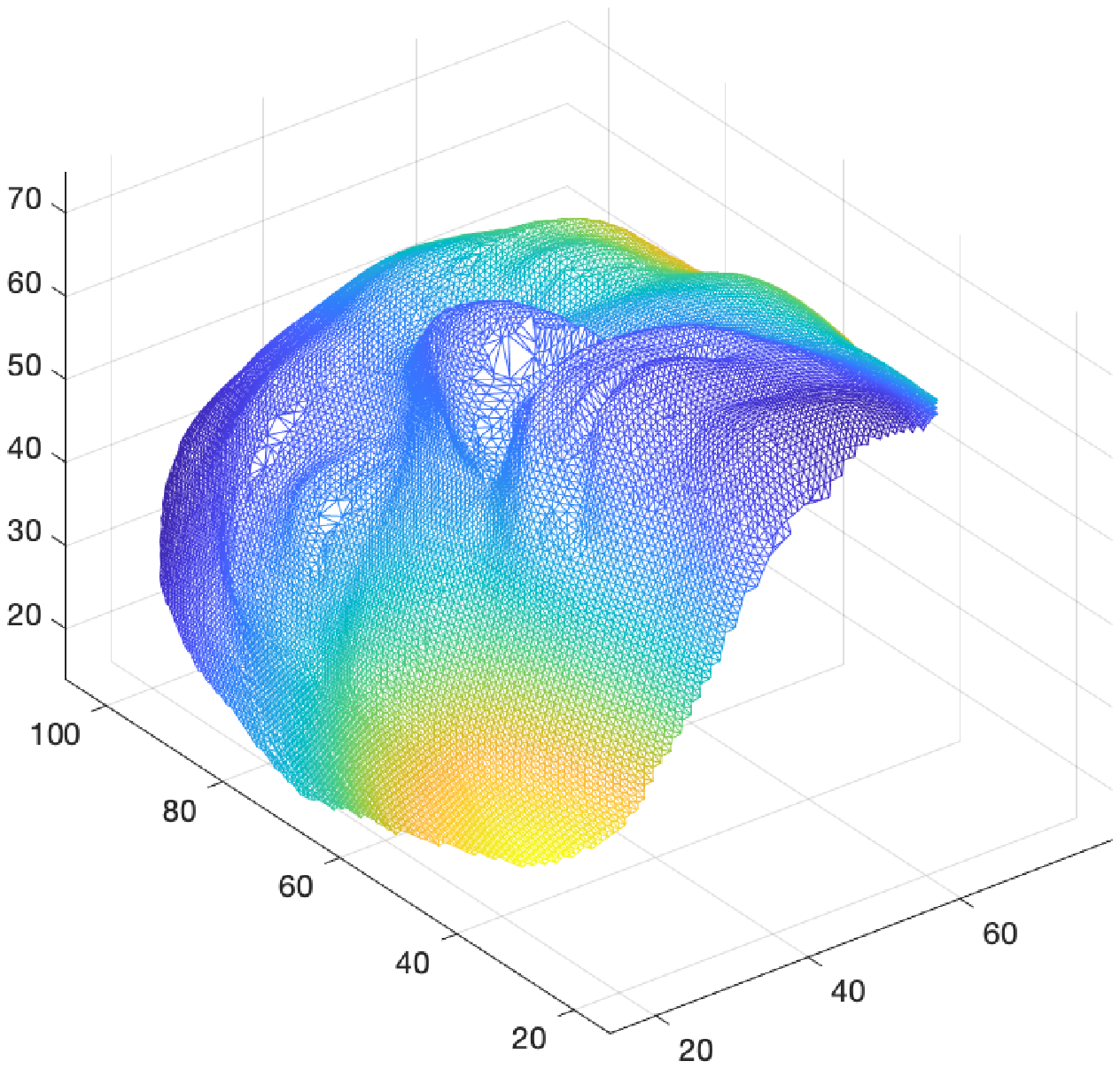}\hskip -12mm
\vskip -10mm
\includegraphics[width=55mm]{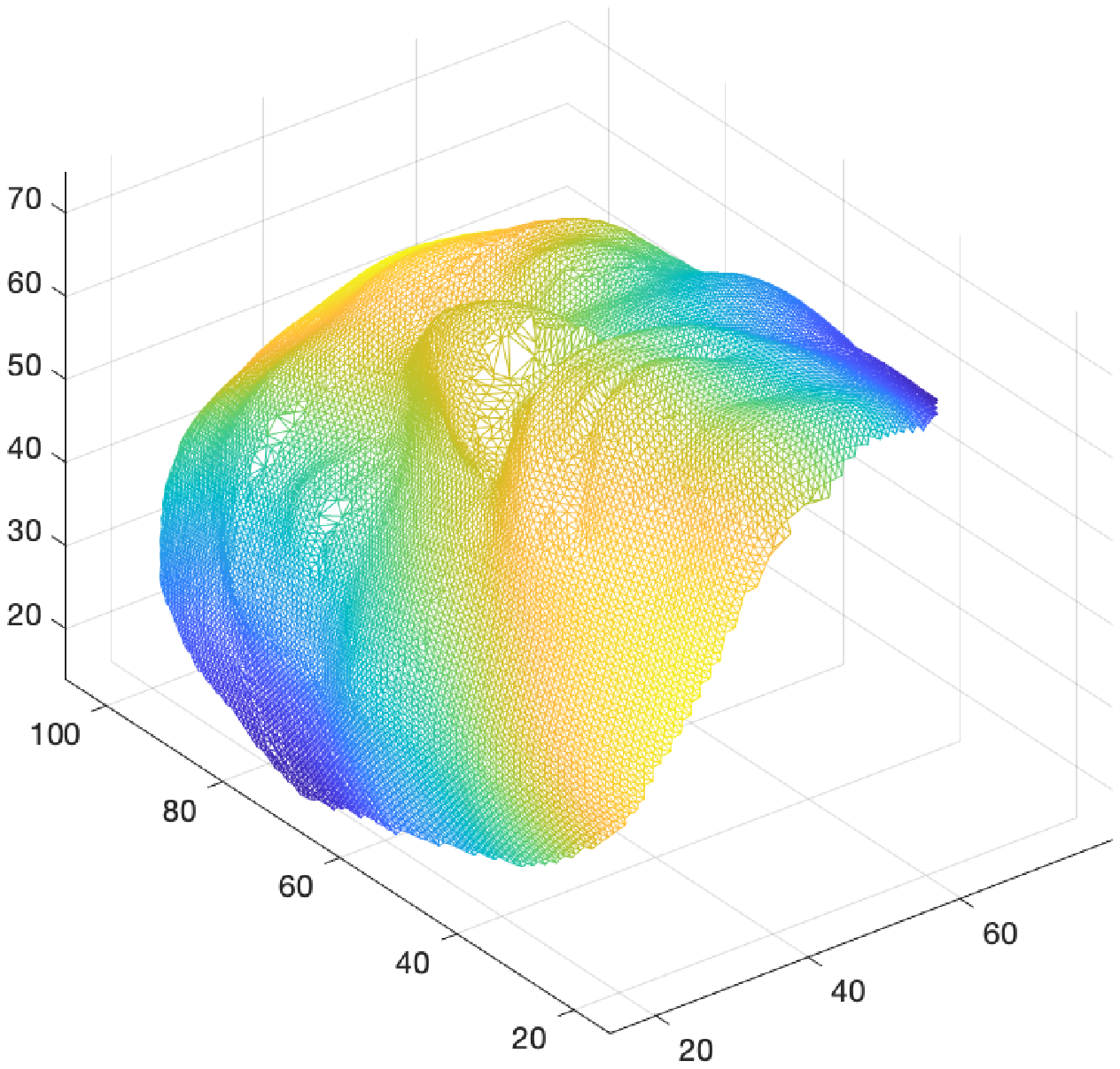}\hskip -12mm
\includegraphics[width=55mm]{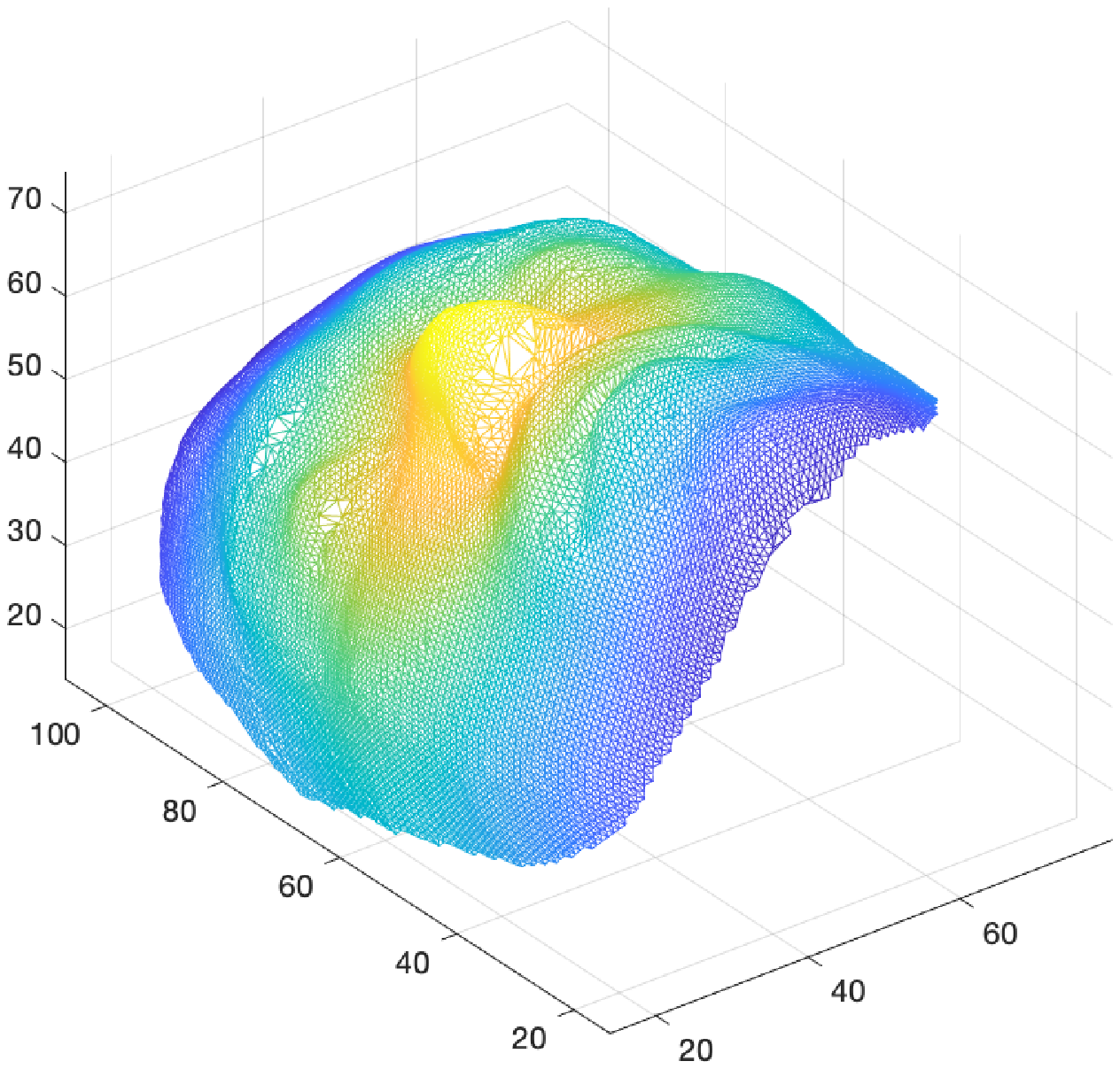}\hskip -12mm
\includegraphics[width=55mm]{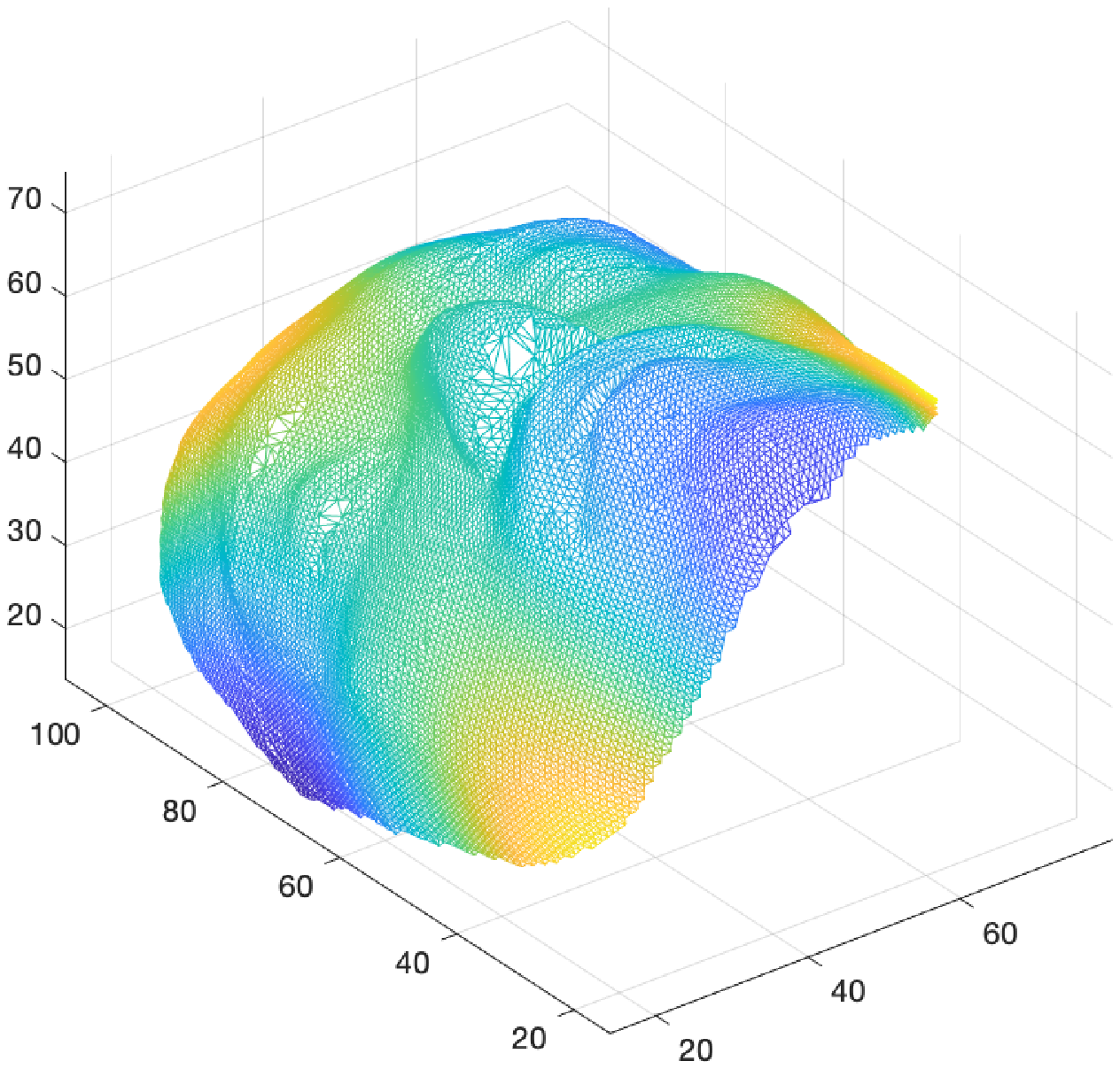}\hskip -12mm
\vskip -10mm
\includegraphics[width=55mm]{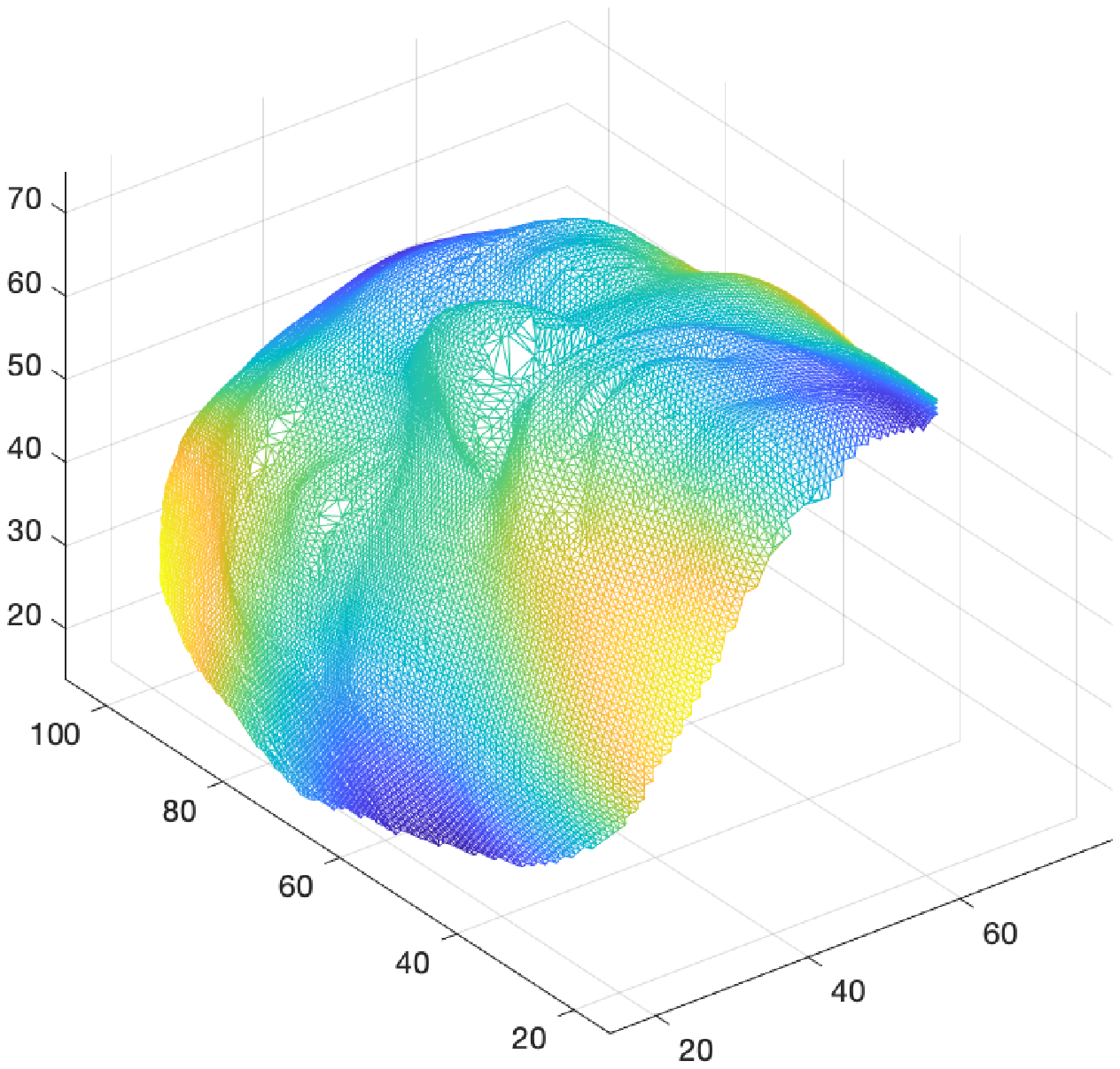}\hskip -12mm
\includegraphics[width=55mm]{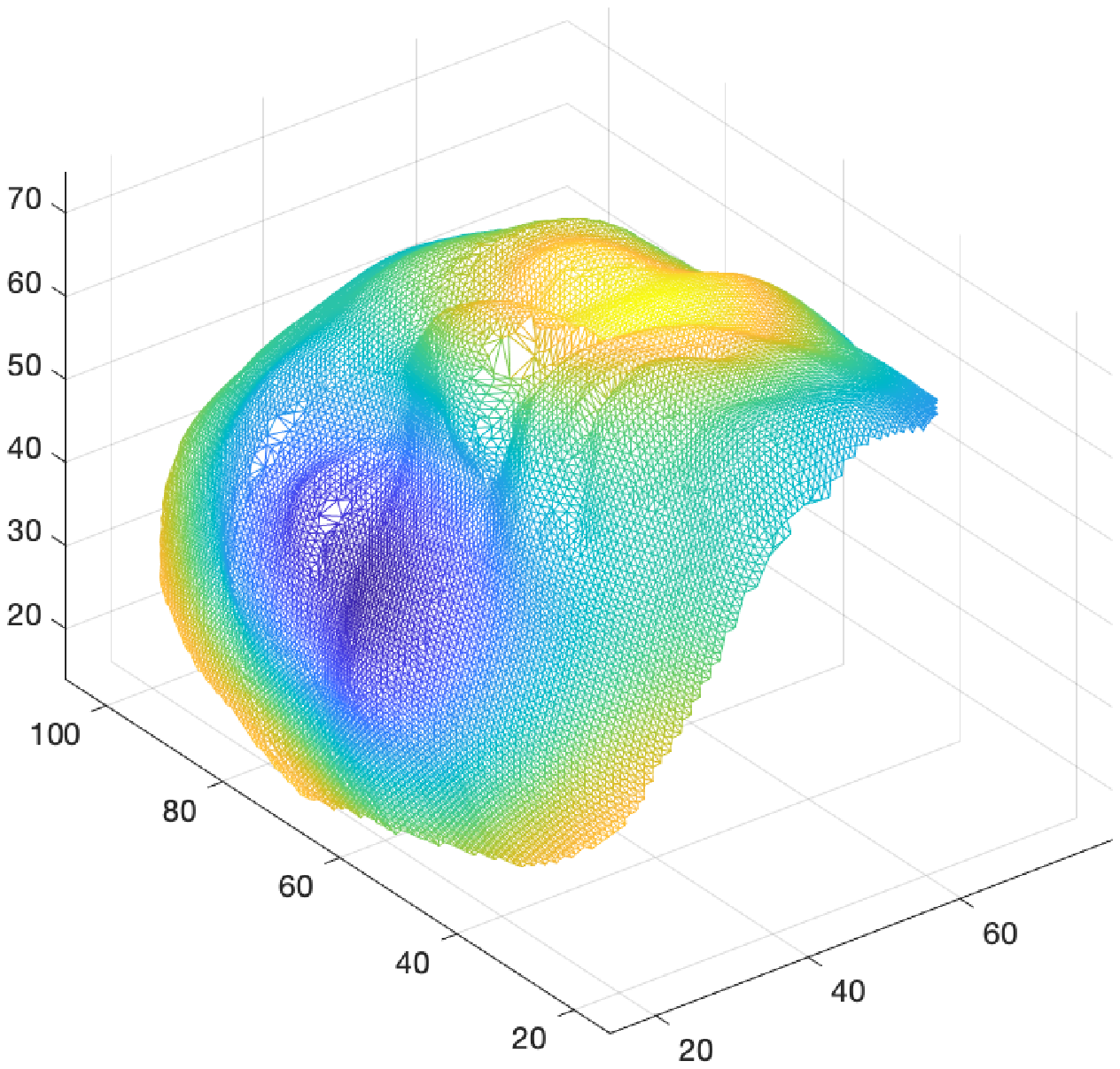}\hskip -12mm
\includegraphics[width=55mm]{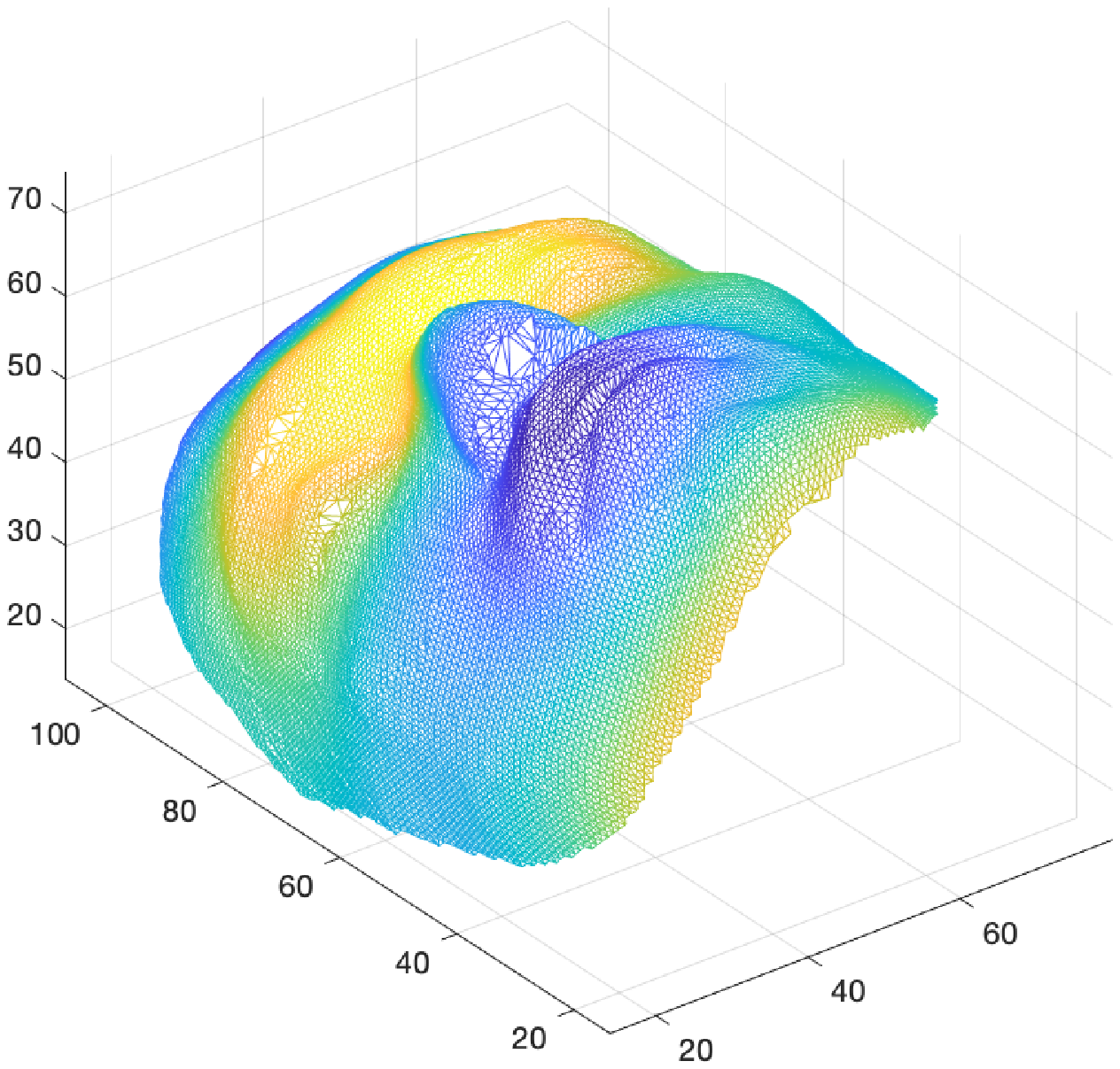}\hskip -12mm
\caption{Eigenfunctions on the face.}\label{fig:face}
\end{figure}

In this section, we report numerical results of the first ten non-zero surface eigenvalues using classical triangle mesh-based algorithm (Algorithm 1), a lattice that is approximated by triangle mesh surfaces (Algorithm 2), and a lattice that is approximated by point cloud (Algorithm 3), at different shell thickness and grid density in Table \ref{tab:face}. We also display values of eigenvectors on the face in Fig.~\ref{fig:face}.

\begin{table}[!htb]
\centering
\centering\caption{The first ten non-zero surface eigenvalues by Algorithm 1, 2, and 3 at different shell thickness and grid density.}\label{tab:tetra}
\begin{tabular}{|l|c|l|l|l|l|l|l|}
\hline
&\multirow{3}*{Surface} &\multicolumn{3}{c|}{Algorithm 2}&\multicolumn{3}{c|}{Algorithm 3}\\
\cline{3-5}\cline{6-8}
&  & $r$ = 0.2 & $r$ = 0.2 &  $r$ = 1
& $r$ = 20 &  $r$ = 20 & $r$ = 30
\\
& & $n$ = 3 & $n$ = 5  & $n$ = 3 & $n$ = 0.05 & $n$ = 0.5 & $n$ = 0.5 
\\
\hline
\multirow{10}*{Eigenvalues}
& -0.0003 & -0.0002 & -0.0002 & -0.0003 & -0.0002 & -0.0002 & -0.0002\\
& -0.0004 & -0.0002 & -0.0003 & -0.0003 & -0.0003 & -0.0003 & -0.0003\\
& -0.0009 & -0.0005 & -0.0007 & -0.0008 & -0.0005 & -0.0007 & -0.0006\\
& -0.0009 & -0.0006 & -0.0007 & -0.0008 & -0.0006 & -0.0007 & -0.0006\\
& -0.0011 & -0.0007 & -0.0009 & -0.0011 & -0.0008 & -0.0009 & -0.0008\\
& -0.0017 & -0.0010 & -0.0013 & -0.0015 & -0.0009 & -0.0013 & -0.0011\\
& -0.0018 & -0.0011 & -0.0014 & -0.0017 & -0.0011 & -0.0014 & -0.0013\\
& -0.0023 & -0.0014 & -0.0019 & -0.0022 & -0.0013 & -0.0019 & -0.0016\\
& -0.0026 & -0.0015 & -0.0020 & -0.0024 & -0.0015 & -0.0020 & -0.0017\\
& -0.0027 & -0.0017 & -0.0021 & -0.0025 & -0.0016 & -0.0021 & -0.0018\\
\hline
\hline
\end{tabular}
\label{tab:face}
\end{table}

\section{Comparison with Analytic Solutions}
This section presents a comparison of the analytic solution and numerical solution of the first ten non-zero Laplacian eigenvalues of $S^2$ and a cone.

\subsection{Comparison with Analytic Solution on $S^2$}\label{sec:comparison_s2}
This section compares our meshless algorithm and Reuter's surface algorithm for Laplacian eigenvalues on $S^2$ with the analytic solution. The results are reported in Table \ref{tab:cone}.

\begin{table}[!htb]
\centering
\caption{Comparison of the first ten non-zero surface Laplacian eigenvalues on $S^2$ computed by the meshless algorithm and surface algorithm with the analytic solution.}
\label{tab:cone}
\begin{tabular}{|l|r|r|r|r|r|r|c|}
\hline
\hline
&\multicolumn{5}{c|}{Meshless Algorithm} & \multirow{2}*{Surface} & \multirow{2}*{Analytic}\\
&$n = 70$&$n = 80$& $n = 90$&$n = 100$&$n = 110$&&\\
\hline
\multirow{10}*{$r = 0.08$} &     -1.855325   &   -1.86539   &   -1.873125   &   -1.879448   &   -1.884542 & -1.9984 & -2\\ 
  &   -1.921102   &   -1.931773   &   -1.938983   &   -1.945495   &   -1.950469 & -2.0005& -2\\ 
  &    -1.935166   &   -1.94506   &   -1.953552   &   -1.959711   &   -1.965312 & -2.0010 & -2\\ 
  &    -5.62263   &   -5.667128   &   -5.696293   &   -5.722359   &   -5.74309  & -5.9924& -6\\ 
  &     -5.634138   &   -5.6763   &   -5.709371   &   -5.733887   &   -5.754284  & -5.9944&-6\\ 
  &    -5.734904   &   -5.754076   &   -5.772994   &   -5.788926   &   -5.801214 & -5.9944& -6\\ 
  &    -5.742877   &   -5.765704   &   -5.784355   &   -5.797827   &   -5.81031  & -5.9966&-6\\ 
  &   -5.821731   &   -5.845335   &   -5.862249   &   -5.877002   &   -5.889156  & -5.9977&-6\\ 
  &     -11.302968   &   -11.376339   &   -11.429968   &   -11.473721   &   -11.507684 & -11.9628& -12\\ 
  &    -11.323203   &   -11.395242   &   -11.451686   &   -11.494617   &   -11.529418 & -11.9690&-12 \\ 
\hline
\hline
\end{tabular}
\end{table}

\subsection{Comparison with Analytic Solution on a Cone}\label{sec:comparison_cone}

The analytic solutions of surface Laplacian eigenvalues on the cone surface are reported in Table \ref{tab:bessel_eigenvalue}, and the comparison with numerical values is reported in Table \ref{tab:cone_multiplicity}.

\begin{table}[!htb]
\centering
\centering\caption{Analytic surface Laplacian eigenvalues on a cone with $n$ from 0 to 9.}\label{tab:bessel_eigenvalue}
\begin{tabular}{|l|cccccccc|}
\hline
\hline
 & \multicolumn{8}{c|}{Analytic Surface Laplacian Eigenvalues $<10.1^2$ from Bessel function $J_n$}\\
\hline
$n$ = 0 & 0 & 3.85927 & 10.3388 & 20.4572 & 34.3839 & 50.7012 & 71.6951 & 94.8567\\
$n$ = 1 & 0 & 3.03793 & 11.4536 & 21.9914 & 36.2099 & 54.2471 & 74.6358 & 99.7367\\
$n$ = 2 & 0 & 8.58743 & 21.7955 & 36.451 & 54.8071 & 76.9932 & 101.448&\\
$n$ = 3 & 0 & 16.485 & 34.7316 & 53.6064 & 76.1614&&&\\
$n$ = 4 & 0 & 26.6572 & 50.1688 & 73.372 & 100.197&&&\\
$n$ = 5 & 0 & 39.0603 & 68.0441 & 95.684&&&& \\
$n$ = 6 & 0 & 53.665 & 88.3116&&&&& \\
$n$ = 7 & 0 & 70.4497&&&&&&\\
$n$ = 8 & 0 & 89.398&&&&&&\\
$n$ = 9 & 0 & &&&&&&\\
\hline
\end{tabular}
\end{table}

\begin{table}[!htb]
\centering
\centering\caption{Analytic surface Laplacian eigenvalues on a cone  with multiplicity with $n$ from 0 to 9 and numerical surface Laplacian eigenvalues by the standard surface algorithm.}\label{tab:cone_multiplicity}
\begin{tabular}{|l|c|c|c|}
\hline
\hline
 & Numerical Eigenvalues 
 & Analytic Eigenvalues & Bessel $J_n$ \\
\hline
\multirow{30}*{Eigenvalues} & 
0 & 0 & 0  \\
&-3.0376&-3.03793 &1\\
&-3.0377 &-3.03793&1\\
&-3.8591&-3.85927 & 0\\
&-8.5843&-8.58743 & 2\\
&-8.5854&-8.58743 & 2\\
&-10.3370&-10.3388 & 0\\
&-11.4486&-11.4536 & 1\\
&-11.4491&-11.4536 & 1\\
&-16.4751&-16.485 & 3\\
&-16.4761&-16.485 & 3\\
&-20.4440&-20.4572 & 0\\
&-21.7744&-21.7955 & 2\\
&-21.7831&-21.7955 & 2\\
&-21.9765&-21.9914 & 1\\
&-21.9777&-21.9914 & 1\\
&-26.6327&-26.6572& 4\\
&-26.6335&-26.6572& 4\\
&-34.3505&-34.3839& 0\\
&-34.6892&-34.7316& 3\\
&-34.6906&-34.7316& 3\\
&-36.1593&-36.2099 & 1\\
&-36.1622&-36.2099 & 1\\
&-36.4016&-36.451& 2\\
&-36.4194&-36.451& 2\\
&-39.0084&-39.0603& 5\\
&-39.0105&-39.0603& 5\\
&-50.0815&-50.1688& 4\\
&-50.0843&-50.1688& 4\\
&-50.6146&-50.7012& 0\\
\hline
\hline
\end{tabular}
\end{table}

\section{Supplementary Material}

Table S1 reports the effect of shell thickness and grid density using Algorithm 2. We observe that denser grids (i.e., greater $n$) lead to an eigenvalue closer to the analytic solution. We also observe that thicker shells lead to eigenvalues better converging to the true value. Table S2 reports analogous results for Algorithm 3.

\setcounter{table}{0}
\renewcommand{\thetable}{S\arabic{table}}

\begin{table}[!htb]
\centering
\caption{The first ten non-zero surface eigenvalues by Algorithm \ref{alg:triangle}.}\label{tab:eig_alg_triangle_nr}
\begin{tabular}{|l|r|r|r|r|r|r|r|r|}
\hline
\hline
&$n = 40$&$n = 50$& $n = 60$&$n = 70$&$n = 80$& $n = 90$&$n = 100$&$n = 110$\\
\hline
  &   -1.3898   &   -1.5193   &   -1.6106   &   -1.6687   &   -1.7121   &   -1.7458   &   -1.7721   &   -1.7938   \\ 
  &   -1.3919   &   -1.5239   &   -1.6116   &   -1.6696   &   -1.7142   &   -1.7468   &   -1.7736   &   -1.7943   \\ 
  &   -1.3964   &   -1.5254   &   -1.6162   &   -1.6704   &   -1.7149   &   -1.7473   &   -1.7743   &   -1.7950   \\ 
   \rowcolor{gray!20} \cellcolor{white}  & -3.6430   &   -4.1117   &   -4.457   &   -4.672   &   -4.8446   &   -4.9748   &   -5.0746   &   -5.1621   \\ 
   \rowcolor{gray!20} \cellcolor{white}  &   -3.6478   &   -4.1221   &   -4.4638   &   -4.6746   &   -4.8487   &   -4.9766   &   -5.0796   &   -5.1637   \\ 
   \rowcolor{gray!20} \cellcolor{white}  &    -4.6407   &   -4.936   &   -5.1297   &   -5.2616   &   -5.3604   &   -5.4318   &   -5.4968   &   -5.5394   \\ 
   \rowcolor{gray!20}  \cellcolor{white} &    -4.646   &   -4.9382   &   -5.1336   &   -5.2632   &   -5.3617   &   -5.4357   &   -5.4972   &   -5.5431   \\ 
   \rowcolor{gray!20} \cellcolor{white}  &   -4.6554   &   -4.9452   &   -5.1414   &   -5.2713   &   -5.3642   &   -5.4379   &   -5.4988   &   -5.5447   \\ 
  \rowcolor{gray!40} \cellcolor{white}  &  -7.9641   &   -8.7511   &   -9.3728   &   -9.7605   &   -10.0492   &   -10.2665   &   -10.4386   &   -10.5886   \\ 
  \rowcolor{gray!40} \cellcolor{white} \multirow{-10}*{$r = 0.02$} &  -7.9974   &   -8.7645   &   -9.3854   &   -9.7646   &   -10.0651   &   -10.2735   &   -10.4463   &   -10.5923   \\ 
\hline
 &   -1.6109   &   -1.6928   &   -1.7454   &   -1.7835   &   -1.8106   &   -1.8321   &   -1.8501   &   -1.8638   \\ 
  &   -1.613   &   -1.6938   &   -1.7466   &   -1.7847   &   -1.8115   &   -1.8334   &   -1.8508   &   -1.8647   \\ 
  &   -1.6144   &   -1.696   &   -1.7481   &   -1.7855   &   -1.8141   &   -1.8347   &   -1.8515   &   -1.8649   \\ 
   \rowcolor{gray!20} \cellcolor{white}  &    -4.4528   &   -4.7585   &   -4.9719   &   -5.1243   &   -5.2314   &   -5.3178   &   -5.3894   &   -5.4423   \\ 
   \rowcolor{gray!20}  \cellcolor{white} &   -4.4577   &   -4.7624   &   -4.9789   &   -5.1253   &   -5.2358   &   -5.3209   &   -5.3919   &   -5.4439   \\ 
   \rowcolor{gray!20} \cellcolor{white}  &   -5.1364   &   -5.3247   &   -5.4322   &   -5.5158   &   -5.5795   &   -5.6271   &   -5.6656   &   -5.6978   \\ 
   \rowcolor{gray!20} \cellcolor{white}  &    -5.1381   &   -5.3262   &   -5.4361   &   -5.5205   &   -5.5821   &   -5.6278   &   -5.6663   &   -5.6989   \\ 
   \rowcolor{gray!20} \cellcolor{white}  &    -5.1424   &   -5.3299   &   -5.4372   &   -5.523   &   -5.583   &   -5.6316   &   -5.6679   &   -5.7003   \\ 
  \rowcolor{gray!40}  \cellcolor{white} &  -9.3796   &   -9.9188   &   -10.2626   &   -10.5272   &   -10.7054   &   -10.8514   &   -10.9721   &   -11.0632   \\ 
  \rowcolor{gray!40} \multirow{-10}*{$r = 0.03$}    \cellcolor{white} &   -9.3941   &   -9.9223   &   -10.2703   &   -10.5338   &   -10.7124   &   -10.8581   &   -10.9768   &   -11.0664   \\ 
\hline
 &   -1.7103   &   -1.7708   &   -1.8104   &   -1.8386   &   -1.8592   &   -1.875   &   -1.8885   &   -1.8983   \\ 
  &   -1.7129   &   -1.7748   &   -1.8115   &   -1.8393   &   -1.8599   &   -1.8759   &   -1.8887   &   -1.8987   \\ 
  &   -1.718   &   -1.7754   &   -1.8127   &   -1.84   &   -1.8605   &   -1.8765   &   -1.889   &   -1.8995   \\ 
   \rowcolor{gray!20} \cellcolor{white}  &    -4.8447   &   -5.0757   &   -5.2287   &   -5.3407   &   -5.4249   &   -5.4896   &   -5.54   &   -5.5824   \\ 
   \rowcolor{gray!20} \cellcolor{white}  &   -4.8491   &   -5.0809   &   -5.2322   &   -5.3446   &   -5.4287   &   -5.4924   &   -5.541   &   -5.5849   \\ 
   \rowcolor{gray!20} \cellcolor{white}  &   -5.356   &   -5.4945   &   -5.5798   &   -5.6394   &   -5.6848   &   -5.7214   &   -5.7524   &   -5.7731   \\ 
   \rowcolor{gray!20} \cellcolor{white}  &   -5.3578   &   -5.4987   &   -5.58   &   -5.6425   &   -5.6867   &   -5.7224   &   -5.7533   &   -5.7749   \\ 
   \rowcolor{gray!20} \cellcolor{white}  &    -5.3669   &   -5.499   &   -5.5827   &   -5.6437   &   -5.6891   &   -5.724   &   -5.7546   &   -5.775   \\ 
   \rowcolor{gray!40} \cellcolor{white}  &   -10.0345   &   -10.4651   &   -10.6959   &   -10.8932   &   -11.0322   &   -11.142   &   -11.2248   &   -11.2982   \\ 
   \rowcolor{gray!40} \multirow{-10}*{$r = 0.04$}   \cellcolor{white} &  -10.0471   &   -10.4685   &   -10.7004   &   -10.8968   &   -11.0383   &   -11.1465   &   -11.2298   &   -11.2995   \\ 
\hline
&   -1.7746   &   -1.8192   &   -1.8496   &   -1.8721   &   -1.8879   &   -1.8997   &   -1.9102   &   -1.9171 \\ 
  &   -1.7756   &   -1.8201   &   -1.8507   &   -1.8723   &   -1.8883   &   -1.9012   &   -1.9109   &   -1.9171   \\ 
  &   -1.7757   &   -1.8208   &   -1.8509   &   -1.8729   &   -1.8887   &   -1.9018   &   -1.9111   &   -1.9171  \\ 
   \rowcolor{gray!20}  \cellcolor{white} &  -5.0739   &   -5.2616   &   -5.3801   &   -5.4703   &   -5.5379   &   -5.5912   &   -5.6303   &   -5.6629 \\ 
   \rowcolor{gray!20}  \cellcolor{white} &   -5.0762   &   -5.2638   &   -5.3831   &   -5.4734   &   -5.5393   &   -5.5923   &   -5.6323   &   -5.6629  \\ 
   \rowcolor{gray!20} \cellcolor{white}  &   -5.5049   &   -5.5981   &   -5.6691   &   -5.7163   &   -5.7516   &   -5.7772   &   -5.8006   &   -5.8117  \\ 
   \rowcolor{gray!20} \cellcolor{white}  &   -5.5091   &   -5.6009   &   -5.6705   &   -5.7197   &   -5.7531   &   -5.7784   &   -5.8008   &   -5.8117   \\ 
   \rowcolor{gray!20}  \cellcolor{white} &   -5.5113   &   -5.6028   &   -5.6724   &   -5.7206   &   -5.7532   &   -5.7819   &   -5.8024   &   -5.8117 \\ 
   \rowcolor{gray!40} \cellcolor{white}  &  -10.4649   &   -10.7581   &   -10.9514   &   -11.1185   &   -11.2235   &   -11.3155   &   -11.3823   &   -11.4291   \\ 
 \rowcolor{gray!40}   \multirow{-10}*{$r = 0.05$}  \cellcolor{white} & -10.4673   &   -10.7602   &   -10.9582   &   -11.1215   &   -11.2279   &   -11.3193   &   -11.3841   &   -11.4291  \\ 
\hline
\hline
\end{tabular}
\end{table}

\begin{table}[!htb]
\centering
\caption{The first ten non-zero surface eigenvalues by Algorithm \ref{alg:meshless}.}\label{tab:eig_alg_meshless_nr}
\begin{tabular}{|l|r|r|r|r|r|r|r|r|}
\hline
\hline
&$n = 40$&$n = 50$& $n = 60$&$n = 70$&$n = 80$& $n = 90$&$n = 100$&$n = 110$\\
\hline
& -1.5006   &   -1.5406   &   -1.566   &   -1.5855   &   -1.5964   &   -1.6086   &   -1.6183   &   -1.6245   \\ 
&-1.7041   &   -1.7527   &   -1.7678   &   -1.796   &   -1.813   &   -1.8199   &   -1.8299   &   -1.8378   \\ 
&-1.7799   &   -1.8225   &   -1.8538   &   -1.8717   &   -1.8915   &   -1.904   &   -1.9131   &   -1.9208   \\ 
\rowcolor{gray!20} \cellcolor{white}& -4.6797   &   -4.8379   &   -4.9042   &   -4.9942   &   -5.0529   &   -5.0702   &   -5.1122   &   -5.1206   \\ 
\rowcolor{gray!20} \cellcolor{white}& -4.8138   &   -5.007   &   -5.1032   &   -5.2002   &   -5.2459   &   -5.2876   &   -5.3188   &   -5.3379   \\ 
\rowcolor{gray!20}\cellcolor{white} & -5.0147   &   -5.1173   &   -5.1679   &   -5.2225   &   -5.2658   &   -5.3115   &   -5.3381   &   -5.373   \\ 
\rowcolor{gray!20} \cellcolor{white}& -5.108   &   -5.1869   &   -5.2726   &   -5.3071   &   -5.3454   &   -5.3699   &   -5.3941   &   -5.4153   \\ 
\rowcolor{gray!20}\cellcolor{white} & -5.2912   &   -5.3956   &   -5.4406   &   -5.4844   &   -5.5344   &   -5.5607   &   -5.5778   &   -5.6016   \\ 
\rowcolor{gray!40}\cellcolor{white} & -9.4901   &   -9.7345   &   -9.7695   &   -9.9365   &   -10.0651   &   -10.0491   &   -10.1249   &   -10.1367   \\ 
\rowcolor{gray!40}\multirow{-10}*{$r = 0.05$}\cellcolor{white} & -9.6765   &   -9.9786   &   -10.1256   &   -10.2391   &   -10.3352   &   -10.407   &   -10.4557   &   -10.4985   \\ 
\hline 
&   -1.6725   &   -1.7074   &   -1.7334   &   -1.7496   &   -1.7626   &   -1.7728   &   -1.7807   &   -1.78681   \\ 
  &   -1.7944   &   -1.836   &   -1.8611   &   -1.8795   &   -1.8949   &   -1.9041   &   -1.9131   &   -1.9199   \\ 
  &   -1.8314   &   -1.8702   &   -1.8979   &   -1.9163   &   -1.929   &   -1.9404   &   -1.9488   &   -1.9563   \\ 
   \rowcolor{gray!20}\cellcolor{white}&  -5.0865   &   -5.2556   &   -5.3538   &   -5.4236   &   -5.4826   &   -5.5197   &   -5.5507   &   -5.5750  \\ 
   \rowcolor{gray!20}\cellcolor{white} &  -5.1252   &   -5.2767   &   -5.3836   &   -5.4605   &   -5.5112   &   -5.5547   &   -5.587   &   -5.6091  \\ 
   \rowcolor{gray!20}\cellcolor{white} &  -5.3581   &   -5.4403   &   -5.4983   &   -5.5376   &   -5.5692   &   -5.5929   &   -5.6156   &   -5.6355   \\ 
  \rowcolor{gray!20}\cellcolor{white} &   -5.3945   &   -5.4805   &   -5.5429   &   -5.5778   &   -5.6091   &   -5.6338   &   -5.6532   &   -5.6711   \\ 
  \rowcolor{gray!20}\cellcolor{white}&   -5.5367   &   -5.6137   &   -5.678   &   -5.7192   &   -5.7502   &   -5.7737   &   -5.7945   &   -5.8127   \\ 
  \rowcolor{gray!40}\cellcolor{white} &   -10.3279   &   -10.589   &   -10.7767   &   -10.8936   &   -10.9816   &   -11.0511   &   -11.1024   &   -11.1439   \\ 
  \rowcolor{gray!40} \multirow{-10}*{$r = 0.06$}\cellcolor{white} & -10.3811   &   -10.6479   &   -10.8146   &   -10.9429   &   -11.0362   &   -11.1066   &   -11.1658   &   -11.2056   \\
  \hline
&   -1.748184   &   -1.7823   &   -1.8033   &   -1.8190  &   -1.8295   &   -1.839443   &   -1.846358   &   -1.851796   \\ 
  &   -1.837259   &   -1.870554   &   -1.8926   &   -1.9078   &   -1.9199  &   -1.9282  &   -1.9360   &   -1.9418   \\ 
  &   -1.854534   &   -1.890521   &   -1.9123   &   -1.9277   &   -1.9404   &   -1.9492  &   -1.9565   &   -1.9628  \\ 
  \rowcolor{gray!20} \cellcolor{white}&   -5.2646   &   -5.4120   &   -5.4916   &   -5.5586   &   -5.6063  &   -5.6409   &   -5.6708   &   -5.6930  \\ 
   \rowcolor{gray!20}\cellcolor{white} &  -5.2797   &   -5.4197   &   -5.5066   &   -5.5717   &   -5.6160  &   -5.6553   &   -5.6812   &   -5.7052   \\ 
  \rowcolor{gray!20}\cellcolor{white} &   -5.50673  &   -5.5794   &   -5.6328   &   -5.6691   &   -5.6955   &   -5.7163   &   -5.7350   &   -5.7481   \\ 
  \rowcolor{gray!20} \cellcolor{white}&   -5.52381  &   -5.6030   &   -5.6544   &   -5.6849   &   -5.7128   &   -5.7325   &   -5.7498   &   -5.7639  \\ 
  \rowcolor{gray!20} \cellcolor{white}&   -5.6335  &   -5.7060   &   -5.7556   &   -5.7861   &   -5.8155   &   -5.8350   &   -5.8534   &   -5.8673   \\ 
  \rowcolor{gray!40}\cellcolor{white} &   -10.6751  &   -10.9116   &   -11.059   &   -11.1716   &   -11.2464   &   -11.3109   &   -11.3573   &   -11.3983   \\ 
  \rowcolor{gray!40} \multirow{-10}*{$r = 0.07$}\cellcolor{white}& -10.6963   &   -10.951   &   -11.0844   &   -11.1932  &   -11.2686   &   -11.3339   &   -11.3842   &   -11.4229   \\ 
  \hline
 &   -1.7956   &   -1.8226   &   -1.8414   &   -1.8553   &   -1.86539   &   -1.8731   &   -1.8794   &   -1.8845   \\ 
  &   -1.859861   &   -1.8890   &   -1.9066   &   -1.9211   &   -1.9317   &   -1.9389   &   -1.9454   &   -1.9504   \\ 
  &   -1.873071   &   -1.9011   &   -1.9226   &   -1.9351   &   -1.9450   &   -1.9535   &   -1.9597   &   -1.9653   \\ 
  \rowcolor{gray!20}\cellcolor{white} &   -5.3766   &   -5.4920   &   -5.5669   &   -5.622   &   -5.6671   &   -5.6962   &   -5.7223  &   -5.7430   \\ 
   \rowcolor{gray!20} \cellcolor{white}&  -5.3863   &   -5.5041   &   -5.5800   &   -5.6341   &   -5.6763   &   -5.7093   &   -5.7338   &   -5.7542   \\ 
   \rowcolor{gray!20}\cellcolor{white} &  -5.597   &   -5.6549   &   -5.7025   &   -5.7349   &   -5.7540   &   -5.7729   &   -5.7889   &   -5.8012   \\ 
   \rowcolor{gray!20}\cellcolor{white} &  -5.6031   &   -5.666   &   -5.7120   &   -5.7428   &   -5.7657   &   -5.7843   &   -5.7978   &   -5.8103   \\ 
   \rowcolor{gray!20}\cellcolor{white} &  -5.6855   &   -5.748   &   -5.7910   &   -5.8217   &   -5.8453   &   -5.8622   &   -5.8770   &   -5.8891   \\ 
   \rowcolor{gray!40} \cellcolor{white}& -10.8871   &   -11.0775   &   -11.2104   &   -11.3029   &   -11.3763   &   -11.4299   &   -11.4737   &   -11.5076   \\ 
  \rowcolor{gray!40}  \multirow{-10}*{$r = 0.08$} \cellcolor{white} & -10.9073   &   -11.1150   &   -11.2387   &   -11.3232   &   -11.3952   &   -11.4516  &   -11.4946  &   -11.5294   \\ 
\hline
\hline
\end{tabular}
\label{tab1}
\end{table}

\clearpage
\mbox{~}

\vskip 1cm

\section*{Appendix}

\vskip 0.5cm

\subsection*{Triangle Mesh of 3-sided and 4-sided Regular Tetrahedron}

The triangle mesh of 3-sided and 4-sided regular tetrahedron we used for Reuter's classical algorithm is displayed as Fig.~A1.

\begin{figure}[H]
\includegraphics[width=0.38\textwidth]{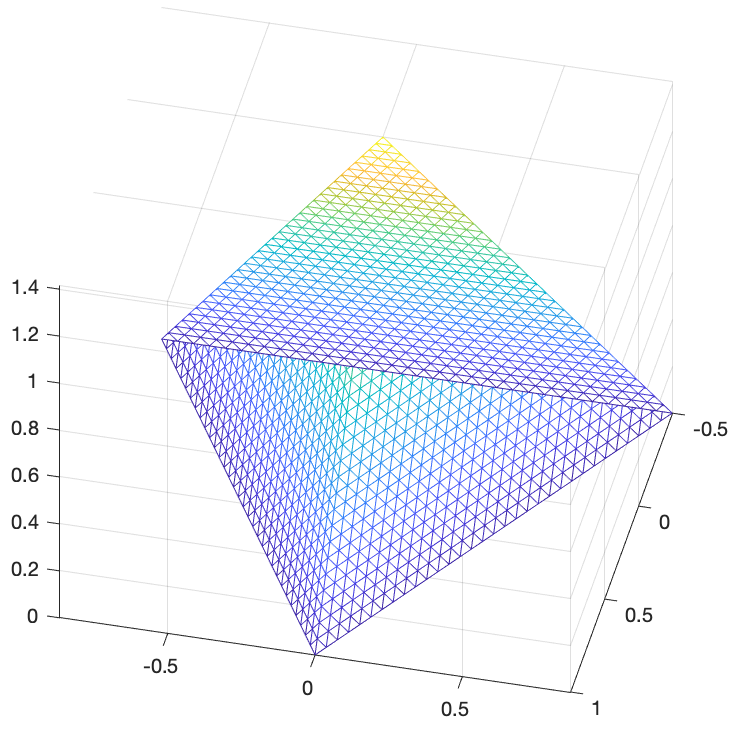}\hskip 10mm
\includegraphics[width=0.38\textwidth]{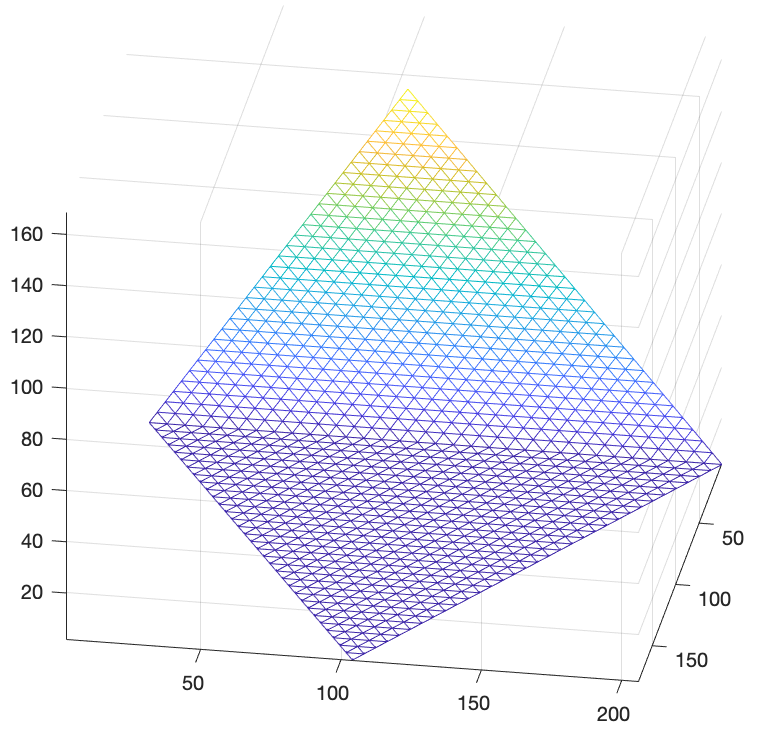}

\noindent \bf Figure A1. \rm Triangle mesh of 3-sided (left) and 4-sided (right) regular tetrahedron.
\end{figure}

\vskip 0.7cm

\subsection*{Cube and Cone}

The point cloud and triangle mesh of a cube and a cone whose vertices are located at the origin and base at $z=1$ with radius = 1.

\begin{figure}[H]
\includegraphics[width=0.35\textwidth]{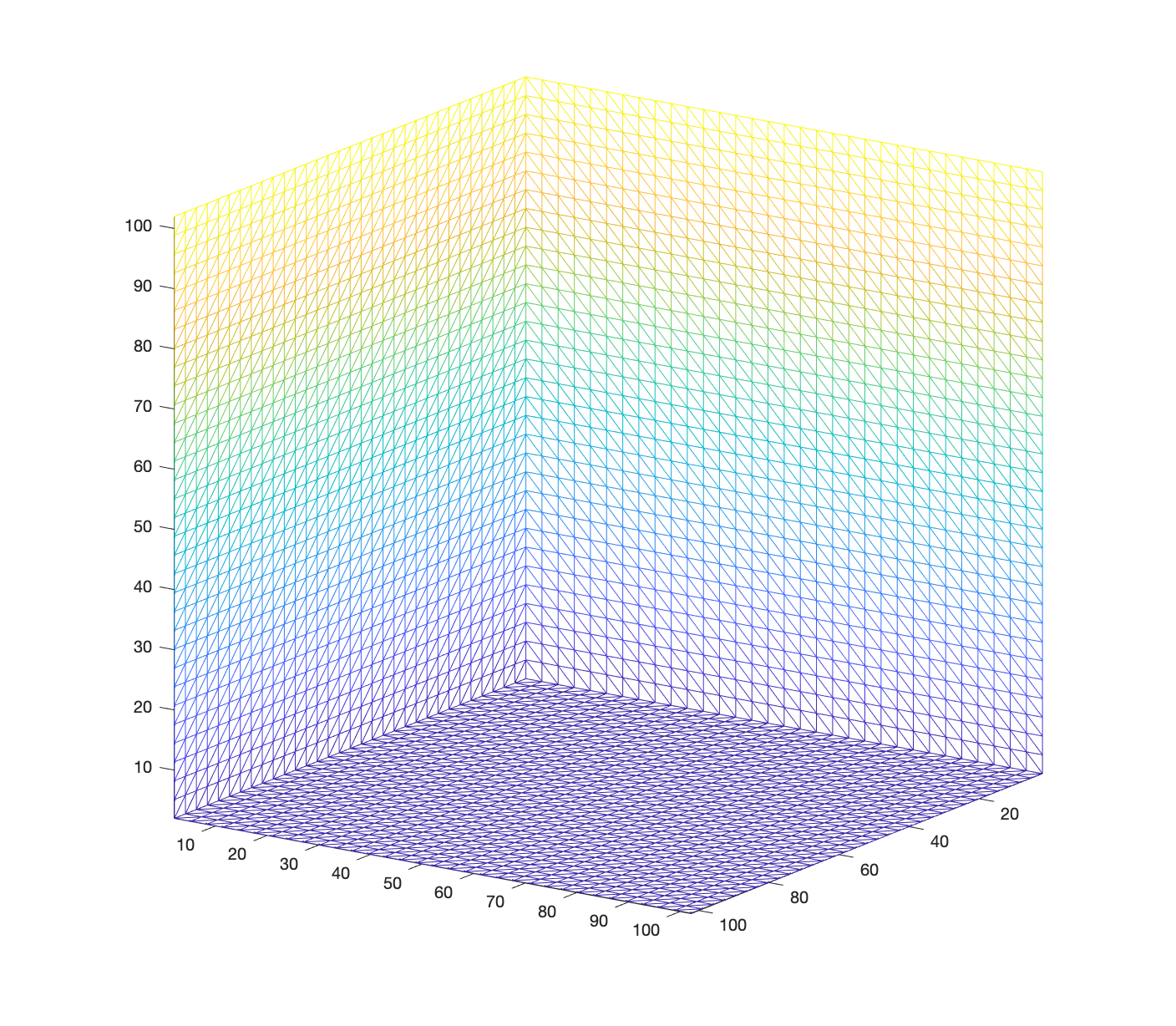}\hskip 10mm
\includegraphics[width=0.35\textwidth]{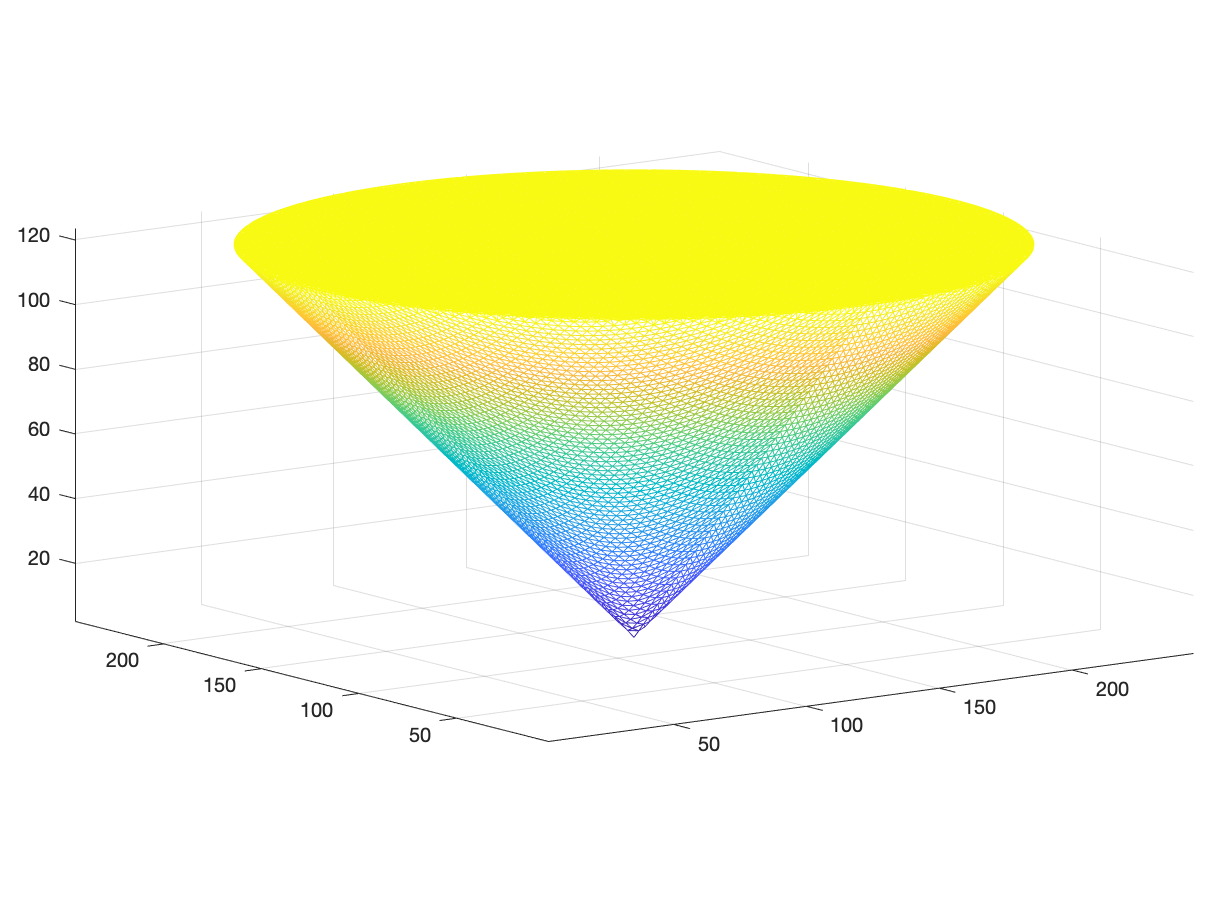}

\noindent \bf Figure A2. \rm Triangle mesh of cube (left) and cone (right).
\end{figure}

\begin{figure}[H]
\includegraphics[width=\textwidth]{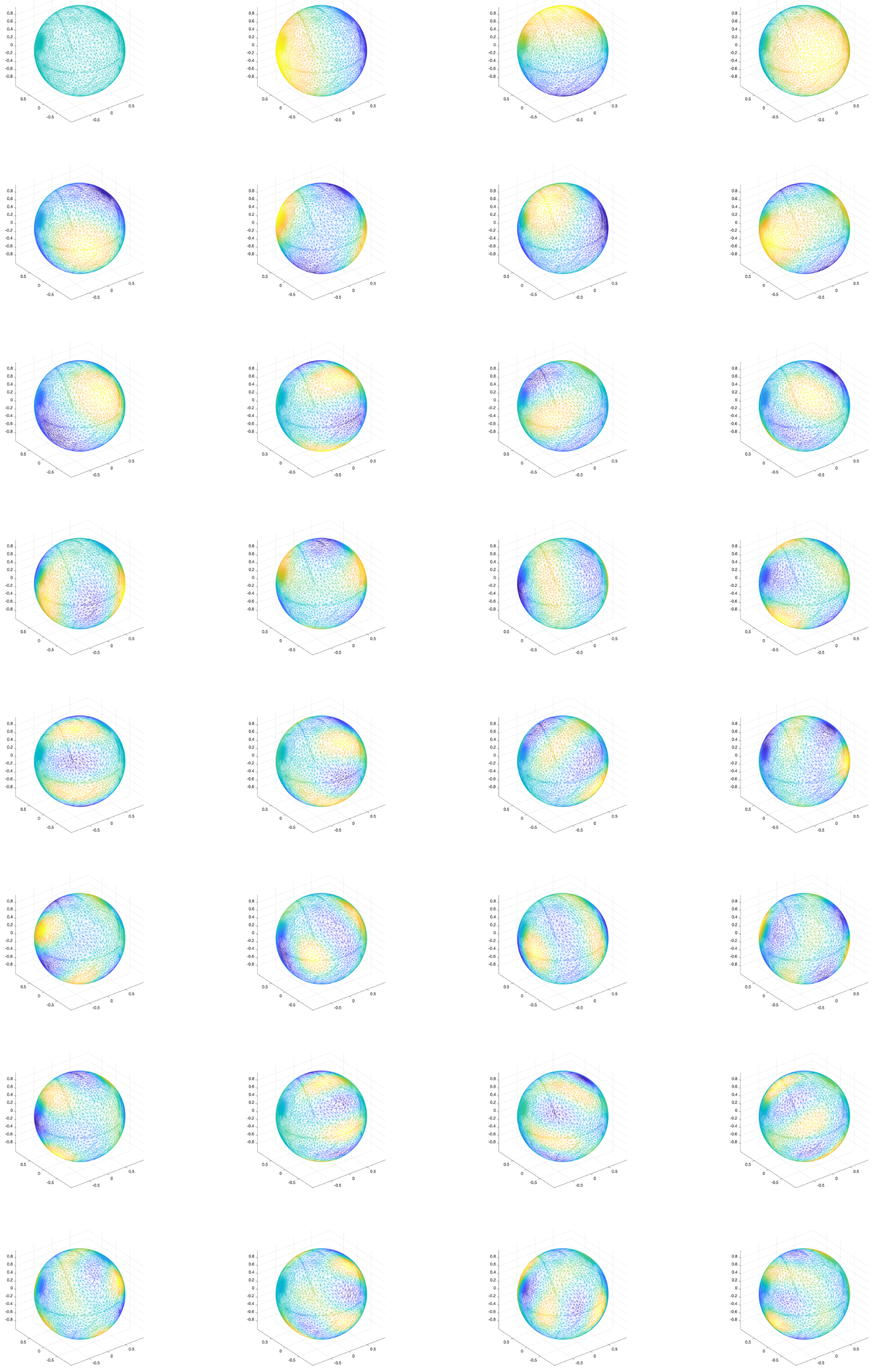}
\noindent \bf Figure A3. \rm Eigenfunctions correspond to first 32 eigenvalues 
on a sphere. \label{fig:func_sphere}
\end{figure}

\begin{figure}[H]
\includegraphics[width=\textwidth]{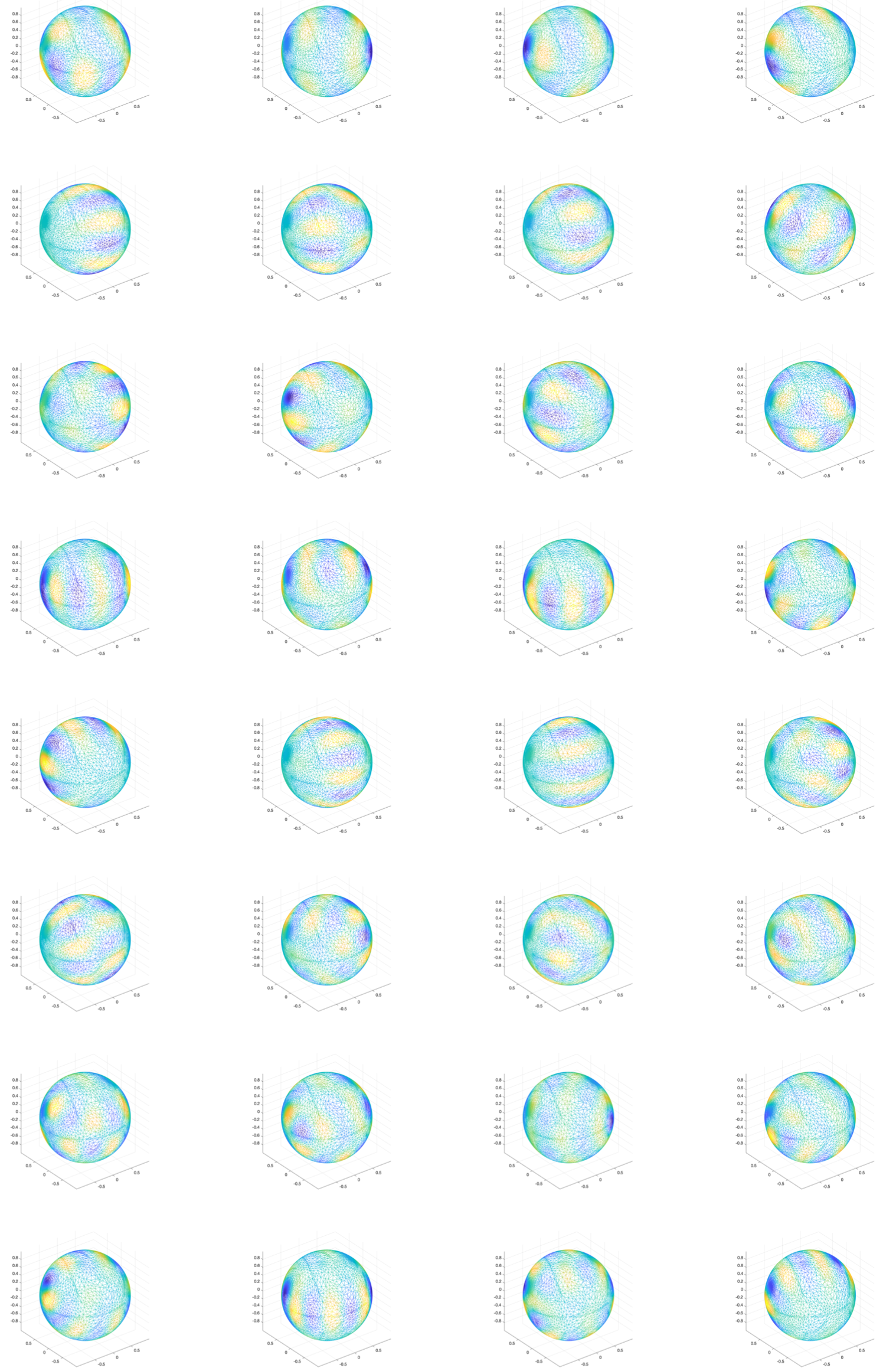}
\noindent \bf Figure A4. \rm Eigenfunctions correspond to the 33rd-64th eigenvalues on a sphere.\label{fig:func_sphere_2}
\end{figure}

\begin{figure}[H]

\includegraphics[width=\textwidth]{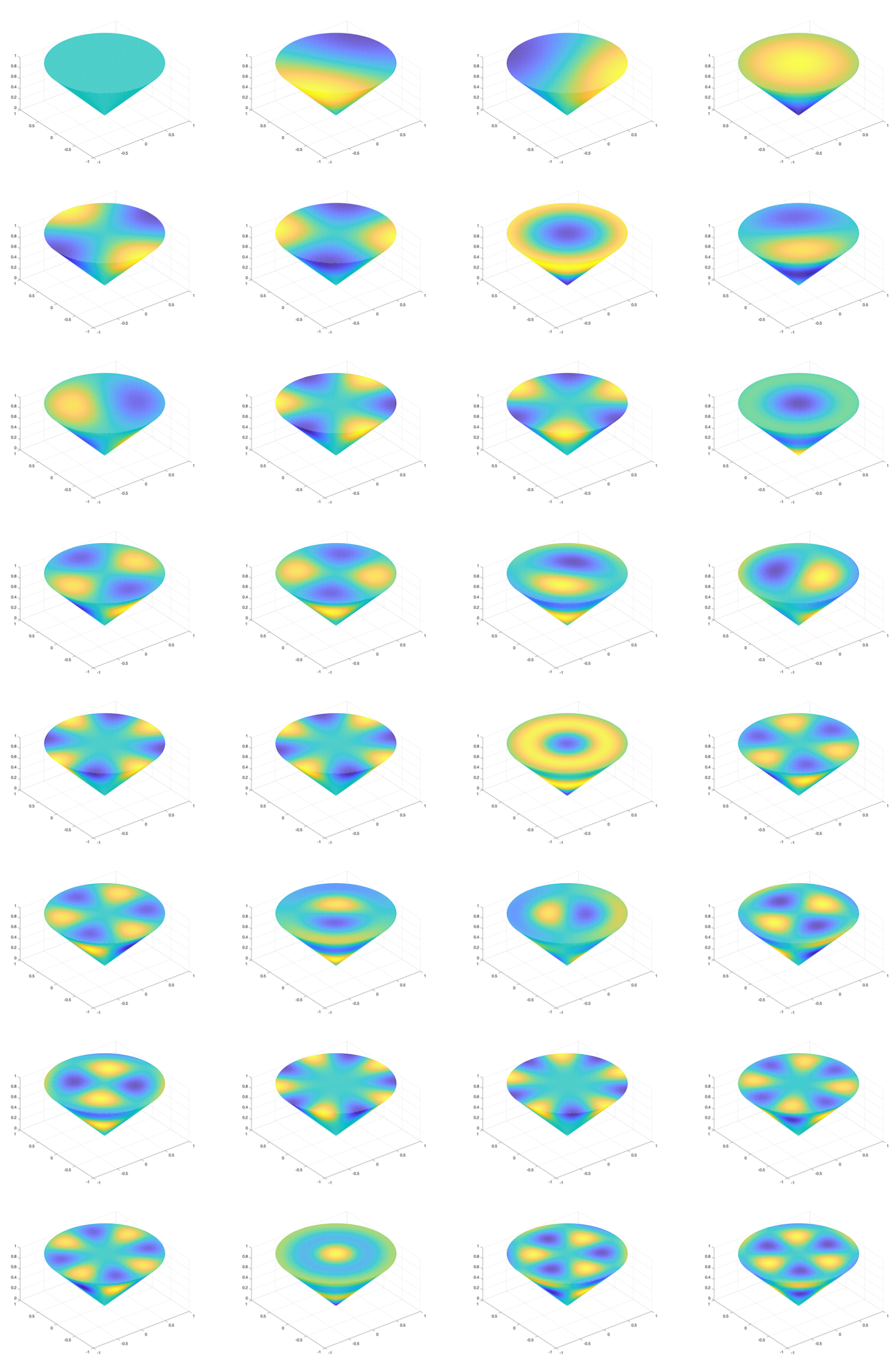}

\noindent \bf Figure A5. \rm Eigenfunctions correspond to first 32 eigenvalues
on a cone.\label{fig:func_cone}
\end{figure}

\begin{figure}[H]
\includegraphics[width=\textwidth]{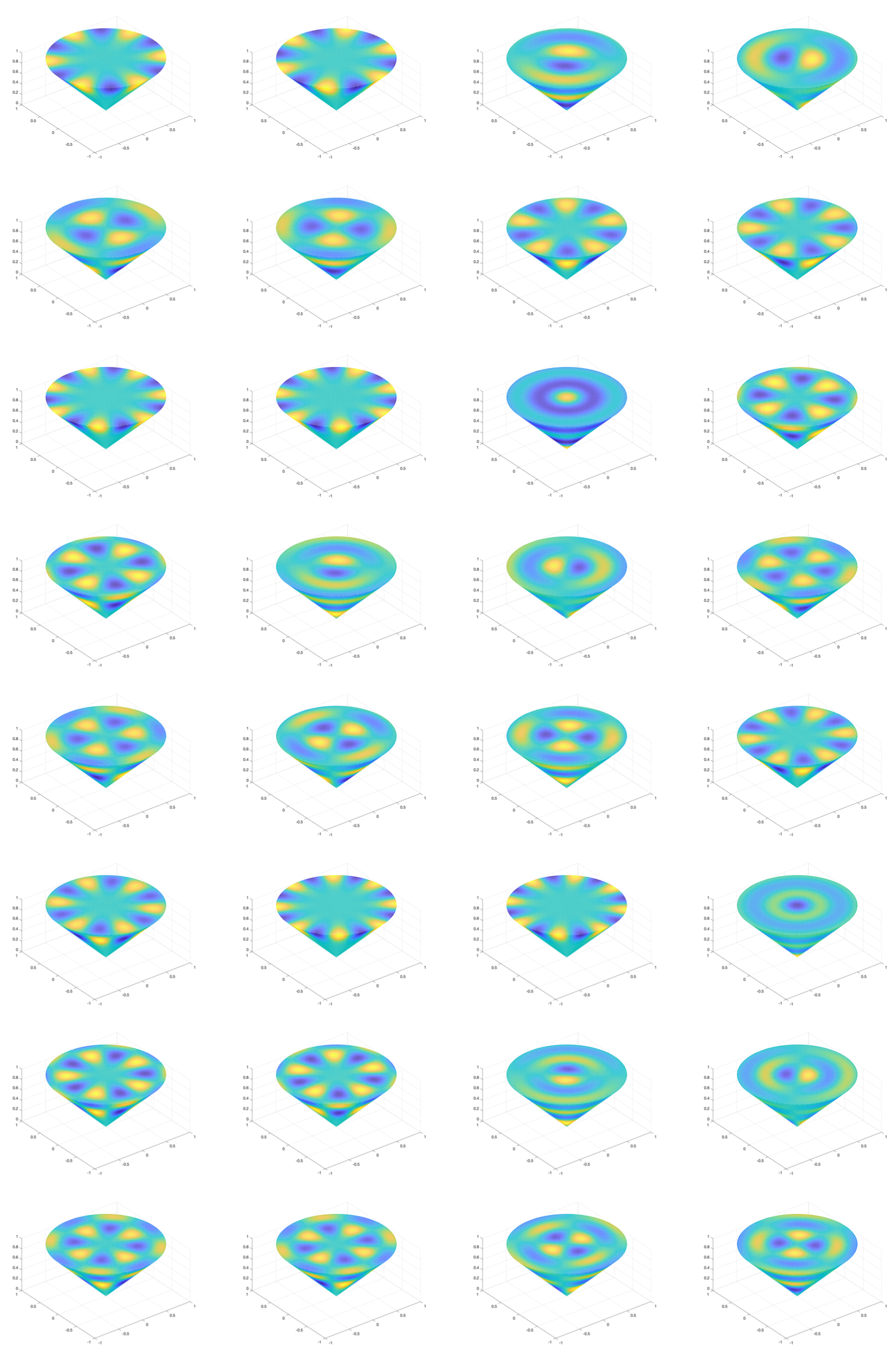}
\noindent \bf Figure A6. \rm Eigenfunctions correspond to the 33rd-64th eigenvalues on a cone.\label{fig:func_cone_2}
\end{figure}

\bibliography{library}
\bibliographystyle{abbrv}

\end{document}